\documentclass[10pt,a4paper]{amsart}
\usepackage{amssymb}
\usepackage{amscd}
\input xypic

\newtheorem{Thm}{Theorem}[section]
\newtheorem{Lem}[Thm]{Lemma}

\newtheorem{Prop}[Thm]{Proposition}
\newtheorem{Defn}[Thm]{Definition}

\def\Alg{Alg_{k}}

\def\N{\mathcal N}
\def\M{\mathcal M}

\def\X{\mathcal X}

\def\R{\mathcal R}

\def\V{\mathcal V}

\def\ext#1#2#3#4{{\rm Ext}^{#1}_{#2}(#3,#4)}

\def\Simp{Simp}

\def\2matr#1#2#3#4{\left(
      \begin{array}{cc}
	#1&#2\\
	#3&#4\\
      \end{array}
      \right)}

\def\matrV#1#2{\left(
      \begin{array}{cc}
	#1\\
	#2\\
      \end{array}
      \right)}
\def\matrH#1#2{\left(
      \begin{array}{cc}
	#1&#2\\
      \end{array}
      \right)}

\def\limprojder#1#2#3
{\mathop{{\lim_{\leftarrow}}^{(#1)}}
_{#2}#3}

\def\limproj#1#2
{\mathop{{\lim_{\leftarrow}}}
_{#1}#2}

\begin{document}

\title{Noncommutative plane curves}
\author{S\o{}ren J\o{}ndrup, Olav Arnfinn Laudal, Arne B. Sletsj\o{}e}

\address{
S\o{}ren J\o{}ndrup\newline
Matematisk institutt\newline
K\o{}benhavns Universitet\newline
Universitetsparken 5\newline
DK-2100 K\o{}benhavn, Denmark}

\address{
Olav Arnfinn Laudal\newline
Matematisk institutt,\newline
University of Oslo,\newline
Pb. 1053, Blindern,\newline
N-0316 Oslo, Norway}

\address{
Arne B. Sletsj\o{}e\newline
Matematisk institutt,\newline
University of Oslo,\newline
Pb. 1053, Blindern,\newline
N-0316 Oslo, Norway}

\email{jondrup@math.ku.dk, arnfinnl@math.uio.no, arnebs@math.uio.no}

\thanks{
This work was partly done during the authors stay at 
Institut Mittag-Leffler.}

\begin{abstract}
    In this paper we study noncommutative plane curves, i.e. 
    non-commutative $k$-algebras for which the 1-dimensional simple 
    modules form a plane curve. We study extensions of simple modules 
    and we try to enlighten the completion problem, i.e. understanding the 
    connection between simple modules of different dimension.
\end{abstract}

\maketitle



\section{Introduction}\label{introduction}
The study of affine algebraic plane curves  goes back to  Descartes and
Fermat, and is still far from complete. There are  many unanswered
natural questions,
one of which is the classification of singularities of such curves. Since
an (embedded) affine algebraic plane curve is the obvious equivalence
class of polynomials
in two variables it is moreover clear that the complete classification of all
isomorphism classes of such curves is an unrealistic task.
In fact, both problems
lead to extremely complex invariant problems, certainly involving
noncommutative
algebraic geometry, see \cite{LP}.
\par
In this paper we shall restrict our attention to a much more
abordable set of problems
regarding affine algebraic plane curves, equally involving
noncommutative algebraic
geometry.
\par
Any affine algebraic plane curve defined on an algebraically closed field $k$,
with coordinate $k$-algebra $A_0$, has many noncommutative models, i.e.
$k$-algebras $A$ such that  $A/([A,A])=A_{0}$. For every such
model A we may consider not only the classical variety
$Simp_1(A)\subset Spec(A_0)$ of
1-dimensional simple representations (the scheme-theoretical closed
points), but also the
$k$-schemes $Simp_n(A)$ of $n$-dimensional simple left representations,
see e.g. \cite{La4}, and the relations between $Simp_{n}(A)$ and $Simp_{m}(A)$.
\par
The structure of these new schemes, naturally associated to affine
algebraic plane curves,
seems to be of interest, not only to noncommutative geometry, but
even to classical
algebraic geometry, see \cite {La3.1} or \cite{La4}.
\par
The purpose of this paper is to start a systematic study of these
noncommutative affine
models of plane quadrics and cubics, and to relate this study to the
relevant literature
on noncommutative algebra, see \cite{Ar}.
\par
We work over an algebraically closed field $k$ of 
characteristic zero, and all $k$-algebras are assumed to be finitely 
generated.


\section{Links to classical ring theory}
In this section we consider simple modules over a PI 
algebra $A$ and prove some results about extensions of 
non-isomorphic simple modules of finite dimension over $k$, for the 
purpose of studying finitely generated indecomposable 
$A$-modules.

\subsection{Extensions of semiprime PI algebras}

A semiprime PI algebra is a $k$-algebra with no nilpotent ideals, 
and satisfying a polynomial identity. Thus the notion generalises the 
concept of an algebraic variety, being defined by a reduced algebra 
satisfying the commutation identity.
\par
For a semiprime 
PI-algebra $A$, the degree $d=deg(A)$ of $A$ is the least integer 
such that $A$ can be embedded in 
a ring of $d\times d$-matrices over a commutative $k$-algebra $R$. 
The degree $d$ is also the maximum of $dim_{k}V$, where $V$ runs 
through all simple $A$-modules (cf. \cite{JJ}), and finally, $d$ is one 
half of the least 
possible degree of a non-zero polynomial $f\in k\langle 
x_{1},\ldots,x_{m}\rangle$ with the property that 
$f(a_{1},\ldots,a_{m})=0$ for all 
$a_{1},\ldots,a_{m}\in A$ and for some $m\ge 1$.
\par
For a simple $A$-module $V$ the left annihilator of $V$
\begin{displaymath}
    \mathfrak{m}=Ann_{A}(V)=\{a\in A\,\vert\, aV=0\}
\end{displaymath}
is a maximal two-sided ideal and 
$A/\mathfrak{m}\simeq End_{k}(V))\simeq M_{l}(k)$, where $l=dim_{k}V$. 
\par
Recall the definition of a $n$-central polynomial.
\begin{Defn}\label{central}
    A noncommutative polynomial $f\in
    k\langle x_{1},\dots,x_{m}\rangle $ is said to be $n$-central if whenever
    evaluated in $M_{n}(R)$ for any commutative $k$-algebra $R$, it is 
    non-vanishing on $M_{n}(R)$ and
    yields values in the center $R$ of $M_{n}(R)$.
\end{Defn}
Let $A$ be a semiprime PI algebra of degree $d=deg(A)$ and $c_{d}$ be a 
$d$-central polynomial. 
Then (cf. \cite{MR}) for all $a_{1},\ldots,a_{m}\in A$ the element 
$c=c_{d}(a_{1},\ldots,a_{m})$ is central in $A$ and non-zero for 
suitable choices of $a_{1},\ldots,a_{m}$. Therefore, if $V$ is a 
simple $A$-module with $dim_{k}V=d$, then for $v\in V$, 
$cv=\overline{c}v=\lambda v$ for some non-zero scalar $\lambda\in k$. 
Here, $\overline{c}$ denotes the image of $c$ in 
$End_{k}(V)\simeq M_{d}(k)$. Moreover, if $W$ is a simple $A$-module and 
$dim\,W<d$, then $cW=0$ (see Chapter 13 of \cite{MR}).
\par
\begin{Prop}\label{PI}
    Let $A$ be a semiprime PI algebra of degree $d$.
    Let $V$ be a simple $A$-module of maximal $k$-dimension and $M$ a 
    finite length module, not containing $V$ as a composition factor.
    Then $Ext_{A}^1(M,V)=Ext_{A}^1(V,M)=0$.
\end{Prop}
\begin{proof}
    We shall prove the vanishing of $Ext_{A}^1(V,M)$. The proof
    for $Ext_{A}^1(M,V)$ is quite similar.
    \par
    Since $Ext_{A}^1(V,-)$ is half-exact as a functor on $A$-$mod$, the
    proposition will follow if we can prove the vanishing result
    for simple modules $W\not\simeq V$.
    \par
    Consider a central element $c$ defined above. Multiplication by 
    $c$ on $V$ (resp. $W$) induces maps $c^{\ast}$ (resp. 
    $c_{\ast}$) on $Ext^1_{A}(V,W)$. Since $c$ is a
    central element of $A$, $c^{\ast}=c_{\ast}$. Notice that 
    \begin{displaymath}
	c\,Ext_{A}^1(V,W)=c^{\ast}(Ext_{A}^1(V,W))=c_{\ast}(Ext_{A}^1(V,W))
    \end{displaymath}
    Since $c$ act as an isomorphism on $V$ we have
    \begin{displaymath}
	c^{\ast}(Ext_{A}^1(V,W))\simeq Ext_{A}^1(V,W))
    \end{displaymath}
    Suppose $W$ is a simple $A$-module and $dim_{k}\,W<d$. Then 
    $c_{\ast}=0$, and so
    \begin{displaymath}
	Ext_{A}^1(V,W)\simeq c_{{\ast}}(Ext_{A}^1(V,W))=0
    \end{displaymath}
    For $W$ of maximal dimension, the element $c$ acts as an
    isomorphism on $W$. Thus, by localisation in $c$ we obtain the 
    implication
    \begin{displaymath}
	Ext_{A}^1(V,M)\not=0\quad\Rightarrow \quad
	Ext_{A[c^{-1}]}^1(V,M)\not=0
    \end{displaymath}
    Since
    $A[c^{-1}]$ is an Azumaya algebra (see \cite{MR}), hence
    has no non-trivial extensions of non-isomorphic simple modules of
    maximal degree, the result follows.
\end{proof}

\subsection{Second layer link}\label{second layer link}

Recall that for prime ideals $\mathfrak{p}$ and $\mathfrak{q}$ of a $k$-algebra $A$, a
second layer link from
$\mathfrak{q}$ to $\mathfrak{p}$, denoted $\mathfrak{q}\leadsto \mathfrak{p}$,
is defined by the existence of an ideal $I$ such that
\begin{displaymath}
    \mathfrak{q}\mathfrak{p}\subseteq I\varsubsetneq \mathfrak{q}\cap \mathfrak{p}
\end{displaymath}
and $(\mathfrak{q}\cap \mathfrak{p})/I$ is a non-trivial torsionfree left 
$A/\mathfrak{q}$-module and right
$A/\mathfrak{p}$-module. In case $\mathfrak{p}$ and $\mathfrak{q}$ are maximal ideals, 
we either have
$\mathfrak{q}\leadsto \mathfrak{p}$ or $\mathfrak{q}\cap \mathfrak{p}=\mathfrak{q}\cdot \mathfrak{p}$.
\par
If $V$ and $W$ are simple $A$-modules of finite dimension over $k$, the
annihilators  $Ann_{A}(V)=\mathfrak{m}$ and $Ann_{A}(W)=\mathfrak{n}$ are maximal ideals.
\begin{Prop}
    Let $V$ and $W$ be non-isomorphic finite dimensional simple $A$-modules
    with annihilators
    $\mathfrak{m}$, resp. $\mathfrak{n}$. Then $Ext_{A}^1(V,W)\not=0$ implies
    that $\mathfrak{n}\leadsto \mathfrak{m}$.
\end{Prop}
\begin{proof}
    Notice that in case $\mathfrak{n}\not\leadsto \mathfrak{m}$, we have
    $\mathfrak{n}\cap \mathfrak{m}=\mathfrak{n}\cdot \mathfrak{m}$, 
    $A/\mathfrak{m}\simeq V^k$ and
    $A/\mathfrak{n}\simeq W^{l}$ for suitable $k$ and
    $l$. Since $V\not\simeq W$ there is a short-exact
    sequence
    \begin{displaymath}
	0\rightarrow Hom_{A}(A,A/\mathfrak{n})\rightarrow
	Hom_{A}(\mathfrak{m},A/\mathfrak{n})\rightarrow
	Ext^1_{A}(V^k,W^{l})\rightarrow 0
    \end{displaymath}
    Observe that
    \begin{align*}
	Hom_{A}(\mathfrak{m},A/\mathfrak{n})
	&\simeq Hom_{A}(\mathfrak{m}/\mathfrak{n}\mathfrak{m},A/\mathfrak{n})\\
	&=Hom_{A}(\mathfrak{m}/\mathfrak{n}\cap \mathfrak{m},A/\mathfrak{n})\\
	&\simeq Hom_{A}(A/\mathfrak{n},A/\mathfrak{n})\\
	&\simeq Hom_{A}(A,A/\mathfrak{n})
    \end{align*}
    Since all modules have finite $k$-dimension we obtain $Ext_{A}^1(V,W)=0$. 
\end{proof}

\subsection{A non PI Result}\label{Soeren}

Let $A=R[\theta;\alpha]$ be a quasipolynomial extension of a 
$k$-algebra $R$, where $\alpha$ is a $k$-automorphism of $R$. Recall 
that a quasipolynomial algebra is the ``ordinary'' polynomial algebra 
over $R$ generated by $\theta$ with relation $\theta \cdot r = \alpha(r)\cdot 
\theta$ for all $r\in R$. 
\begin{Prop}
    Let $V$ and $W$ be two simple non-isomorphic $A$-modules with 
    $\theta V=\theta W=0$. Let $Ann_{R}(V)=\mathfrak{m}$ and 
    $Ann_{R}(W)=\mathfrak{n}$. 
    Then the following are equivalent:
    \begin{itemize}
	\item[i)] $Ext_{A}^1(V,W)\neq 0$
	\item[ii)] $Ext_{R}^1(V,W)\neq 0$ or $\alpha(\mathfrak{m})=\mathfrak{n}$.
    \end{itemize}
\end{Prop}
\begin{proof}
    It is straightforward to see that for two $A$-modules $V$ and $W$ 
    with $\theta V=\theta W=0$, $Ext_{R}^1(V,W)\neq 0$ implies 
    $Ext_{A}^1(V,W)\neq 0$. Therefore it suffices to show that in case
    $Ext_{R}^1(V,W)= 0$, we have $Ext_{A}^1(V,W)\neq 0$ if and only if 
    $\alpha(\mathfrak{m})=\mathfrak{n}$.
    \par
    Suppose $R/\mathfrak{m}\simeq V^k$ where $k=dim\,V$. By assumption  
    the short-exact sequence
    \begin{displaymath}
	0\rightarrow (\mathfrak{m},\theta)\longrightarrow A\longrightarrow 
	V^k\rightarrow 0
    \end{displaymath}
    induces another short-exact sequence
    \begin{displaymath}
	0\rightarrow Hom_{A}(A,W)\buildrel{\phi}\over\longrightarrow 
	Hom_{A}((\mathfrak{m},\theta),W)
	\longrightarrow (Ext_{A}^1(V,W))^k \rightarrow 0
    \end{displaymath}
    and we observe that $Ext_{A}^1(V,W)=0$ if and only if $\phi$ is onto.
    \par
    Suppose $f\in Hom_{A}(A,W)$, then since $\theta W=0$ we have 
    $f(\theta)=0$, thus 
    \begin{displaymath}
	Ext_{A}^1(V,W)\neq 0
    \end{displaymath}
    if we can find an element $g\in 
    Hom_{A}((\mathfrak{m},\theta),W)$ with $g(\theta)\neq 0$.
    \par
    Since $Ext_{R}^1(V,W)=0$ we have an isomorphism 
    $Hom_{R}(R,W)\simeq Hom_{R}(\mathfrak{m},W)$ and the restriction of 
    $g\in  Hom_{A}((\mathfrak{m},\theta),W)$ to $\mathfrak{m}$ is 
    multiplication by some element $w\in W$. Hence for any $b(\theta)\in 
    A$ and $m\in\mathfrak{m}$
    \begin{align*}
	g(rm+b(\theta)\theta)
	&=rg(m)+b(\theta)g(\theta)\\
	&=rmw=(rm+b(\theta)\theta)w
    \end{align*}
    and $g$ is clearly in the image of $\phi$.
    Consequently 
    $Ext_{A}^1(V,W)\neq 0$ is equivalent to the existence of
    $g\in Hom_{A}((\mathfrak{m},\theta),W)$ with $g(\theta)\neq 0$.
    \par
    For any $g\in Hom_{A}((\mathfrak{m},\theta),W)$ and $m\in\mathfrak{m}$ we have 
    \begin{displaymath}
	0=\theta\cdot g(m)=g(\theta\cdot m)=g(\alpha(m)\cdot \theta)=
	\alpha(m)\cdot g(\theta)
    \end{displaymath}
    If $g(\theta)\neq 0$, then $\alpha(\mathfrak{m})\subset 
    Ann_{R}(g(\theta))=Ann_{R}(W)=\mathfrak{n}$ since $W$ is simple. But 
    then $\alpha(\mathfrak{m})= \mathfrak{n}$, since $\alpha$ is an 
    automorphism of $R$ and $\mathfrak{m}$ and $\mathfrak{n}$ are maximal 
    ideals.
    \par
    On the other hand, if $\alpha(\mathfrak{m})=\mathfrak{n}$, then for any 
    $0\neq w\in W$, we define 
    \begin{displaymath}
	g:(\frak{m},\theta)\rightarrow W
    \end{displaymath}
    by
    \begin{displaymath}
	g(m+(\sum_{j}b_{j}\theta^j)\theta)=b_{0}w
    \end{displaymath}
    where $m\in\frak{m}$ and $b_{j}\in R$.
    This is a well-defined, non-zero $A$-homomorphism. In fact $g$ is 
    $R$-linear and
    \begin{displaymath}
	\theta g(m+(\sum_{j}b_{j}\theta^j)\theta)=\theta b_{0}w=0
    \end{displaymath}
    and
    \begin{align*}
	g(\theta(m+(\sum_{j}b_{j}\theta^j)\theta))
	&=g(\alpha(m)\theta+\theta (\sum_{j}b_{j}\theta^j)\theta)\\
	&=(\alpha(m)+\theta (\sum_{j}b_{j}\theta^j))g(\theta)\\
	&=\alpha(m)g(\theta)+\theta(\sum_{j}b_{j}\theta^j)g(\theta)=0
    \end{align*}
\end{proof}


\section{Geometrical approach}\label{classical}

In this section we explore two approaches to understanding the 
geometry of the scheme $Simp_{n}(A)$, the global approach using 
invariant theory, and the local approach using deformation theory.

\subsection{Representations of algebras. }\label{rep}

Let $k$ be an algebraically closed field and $A$ a finitely
presented associative $k$-algebra.
A left $A$-module structure on a $k$-vector space $M$ is given by a
$k$-algebra homomorphism
\begin{displaymath}
    \phi:A\longrightarrow {End}_{k}(M)
\end{displaymath}
This homomorphism $\phi$ is called a representation, or equivalently,
the structure map of the module $M$. Two structure maps give
isomorphic $A$-module structures on $M$ if they
differ by a ring automorphism of ${End}_{k}(M)$. Such
automorphisms are conjugations by an invertible linear operator
$T:M\longrightarrow M$.
Thus the set of left $A$-module structures
on $M$ is the set of structure maps $\phi:A\longrightarrow
{End}_{k}(M)$, modulo conjugation.
\par
By choosing a $k$-basis for $M$ an n-dimensional representation of $A$ 
is thus given by a
homomorphism of $k$-algebras,
\begin{displaymath}
    \phi: A\longrightarrow M_{n}(k)
\end{displaymath}
where $M_{n}(k)$ is the $k$-algebra of $n\times n$-matrices over $k$.
The representation $\phi$ is simple if it is surjective.
In the finite dimensional case, the automorphism group acting on the
set $\X_{A}^{(n)}:={Mor}_{\Alg}(A,M_{n}(k))$ of $n$-dimensional
representations is
the group $PGl_{n}(k)$. Since $A$ is finitely presented it is easy
to see that
$\X_{A}^{(n)}$ is a well-defined affine variety.
\par
Let $\{x_{i}\}_{i=1,...,m}$ be generators of $A$ as a $k$-algebra. 
Associate to the representation
$\phi\in \X_{A}^{(n)}$, the point of the affine space $\mathbf{
A}^{m\cdot n^2}$ with coordinates
$\{a^i_{p,q}\},\ i=1,\ldots,m,\ p,q=1,\ldots,n,$ given by 
$\phi (x_{i})=(a^i_{p,q})$.
Let $f_1,\ldots,f_{r}$ be the relations defining $A$, and
put $\Gamma =k[x^i_{p,q}]/(\tilde{f}_1,\ldots,\tilde{f}_{r})$, where 
$\tilde{f}_{j}=f_{j}((x_{p,q}^1),\ldots,(x_{p,q}^m))$. Then
$\X_{A}^{(n)}=Spec(\Gamma)$,
and we have a versal family of $n$-dimensional $A$-modules with
basis $\X_{A}^{(n)}$ 
given by the representation,
\begin{displaymath}
    \Phi:A\longrightarrow M_{n}(\Gamma)
\end{displaymath}
defined by,
\begin{displaymath}
    \Phi(x_{i})=(cl(x^i_{p,q}))
\end{displaymath}
where $cl(x^i_{p,q})$ denotes the class of $x^i_{p,q}$ in $\Gamma$.
Clearly $PGl_{n}(k)$ acts on $\X_{A}^{(n)}$ and on the k-algebra $\Gamma$. Put
\begin{displaymath}
    Repr_{n}(A)=\X_{A}^{(n)}/PGl_{n}(k).
\end{displaymath}
The underlying set of this quotient is the set of orbits of
$\X_{A}^{(n)}$ under the
action of $PGl_{n}(k)$, i.e. the set of isomorphism classes of $n$-dimensional
representations of $A$.
The set of isomorphism classes of simple $n$-dimensional
representations, denoted by
$Simp_{n}(A)$ is a subset of $Repr_{n}(A)$.
\par
Notice that the prime spectrum $Spec (A)$ of a ring $A$ is the set of prime
ideals with the Jacobson topology. Clearly $Simp_{n}(A)$ is also a
subset of $Spec(A)$.

\subsection{The free noncommutative algebra in m variables}\label{spacesimp}

Let $S=k\langle x_{1},\dots,x_{m}\rangle$ be the free algebra in $m$
variables. A representation $\phi:S\rightarrow
M_{n}(k)$ of $S$ is given by a set of $n\times n$-matrices
$\phi(x_{1}),\dots,\phi(x_{m})$ and we identify the variety $\X_{S}^{(n)}$
with the
affine space $\textbf{A}^{mn^2}$ of $m$-tuples of $n\times
n$-matrices. In general there is no way of describing
$Repr_{n}(S)$ as an algebraic scheme, but restricting to the simple
modules this is possible, as shown in \cite{Pr0}.
\par
Clearly,  $\Gamma:=\Gamma({\mathcal O}_{\X_{S}^{(n)}})\simeq k[x^i_{s,t}]$
is the free commutative $k$-algebra of polynomials in the
variables $x^i_{s,t}$ for $i=1,\dots,m$ and $s,t=1,\dots,n$.
The versal family in this case is
\begin{displaymath}
    \Phi:S=k\langle x_{1},\dots,x_{m}\rangle \longrightarrow
    M_{n}(\Gamma
    )
\end{displaymath}
\noindent
given  by $x_{i}\mapsto M_{i}$, where $M_{i}$ is a matrix whose
$s,t$-entry is $x^i_{s,t}$. The images
$\pi(S)$ of $\pi$ is the ring of generic $n\times n$-matrices and the
generators $M_{1},\dots,M_{m}$ are in the literature called the
generic matrices over $k$.
\par
The automorphism group $G:=PGl_{n}(k)$ acts on the ring $\Gamma$ by
conjugation. Moreover $G$
acts on $M_{n}(\Gamma)$ by double conjugation, leaving $\Phi$ invariant.
Therefore it is clear that the
coefficients of the characteristic polynomials of the generic
matrices are invariant, i.e. that they are contained in $\Gamma^G$. 
However, there is no versal family defined over $\Gamma^{G}$, 
extending the simple modules, see \cite{La4}.
Define $C_{n}$ to be the subring of $\Gamma\simeq k[x^i_{s,t}]$ generated by
the coefficients of the characteristic polynomials of the generic
matrices. This ring is the trace ring, or the trace ring of dimension $n$
to be more precise.
\par
Artin \cite{Ar} conjectured and Procesi \cite{Pr1} proved that
the trace ring $C_{n}$ is precisely the subring $\Gamma^G$ of $k[x^i_{s,t}]$.
\par
Procesi has also shown \cite{Pr1} that the closed points of $Spec(C_{n})$ are in
1-1 correspondance with equivalence classes of semisimple representations
under the natural map
\begin{displaymath}
    [\rho_{\gamma}]\longmapsto \rho_{\gamma}(\Phi)
\end{displaymath}
where $\rho_{\gamma}:C_{n}\rightarrow k(\gamma):=k$ is a closed point
of $Spec(C_{n})$ and $\rho_{\gamma}(\Phi):S\rightarrow M_{n}(k)$ is 
the corresponding $n$-dimensional representation.
\par
There is a subset $F_{n}(S)\subset S$, see e.g. sect. 13.7 of \cite{MR}, called
the Formanek center of $n$-central polynomials such that
$\phi:S\rightarrow M_{n}(k)$
is a simple representation if and only if $\phi(F_{n}(S))\not=0$. The
set of prime ideals $\mathfrak{p}$ in $S$ such that $\mathfrak{p}\not\supset F_{n}(S)$ is
an open set of $Simp_{n} (S)$. The subset $F_{n}(S)\subset S$ consists of 
$n$-central
elements and by definition $\Phi(F_{n}(S))$ sits in the center of
$M_{n}(\Gamma)$. The center is $PGl_{n}(k)$-invariant and $\Phi(F_{n}(S))$,
being invariant, is a subset $\Phi_{n}$ of the trace ring $C_{n}$. A  
homomorphism
$\psi:C_{n}(S)\rightarrow k$ corresponds to a simple
representation if and only if $\psi(\Phi_{n})\not=0$. Procesi has 
proved the following result:
\begin{Thm}\label{proc}{(Procesi)}
    The variety $Spec ((C_{n})_{\Phi_{n}})$ of simple $n$-dimensional
    representations of $S$ is a smooth variety of dimension $m
    n^2 - (n^2-1)=(m-1)n^2 + 1$
\end{Thm}
Unluckily it is, in general, very difficult to compute $F_{n}(S)$, and
therefore also
$Simp_{n} (S)$, using this method.

\subsection{The free noncommutative algebra in two variables}\label{planesimp}

As an example, let us consider the noncommutative affine plane
$S=k\langle x,y\rangle $. A left $S$-module
$M$, isomorphic to $k^2$ as $k$-vector space, is given by a ring
homomorphism
\begin{displaymath}
    \phi:S\longrightarrow {\rm End}_{k}(M)\simeq M_{2}(k)
\end{displaymath}
\noindent The $S$-module $M$ is simple if and only if $\phi$ is
surjective, i.e. iff $X:=\phi(x)$ and $Y:=\phi(y)$ generate
$M_{2}(k)$ as $k$-algebra. It is easily  checked that this is 
equivalent to $\det [X,Y]\not=0$ (\cite{AS}).
\par
In this case the trace ring
\begin{displaymath}
    C_{2}=k[t_{X}, d_{X}, t_{Y}, d_{Y}, t_{XY}]
\end{displaymath}
is the polynomial ring on the five trace elements.
We use the notation $t_{P}$ for the trace and
$d_{P}$ for the determinant of the matrix $P$. Notice that the
determinant itself is a trace, due to the formula
\begin{displaymath}
    d_{P}=\frac{1}{2}(t_{(P^2)}-(t_{P})^2)
\end{displaymath}
\noindent where we, of course have to assume char$(k)\not= 2$.
\par
Notice that in this case, the Formanek center consists of one single 
element, and non-vanishing of this Formanek element, 
defined by
\begin{displaymath}
    d_{[X,Y]}=-\frac{1}{4}
    ((2t_{XY}-t_{X}t_{Y})^2-(t_{X}^2-4d_{X})(t_{Y}^2-4d_{Y}))
\end{displaymath}
is equivalent to simplicity of the
$S$-module $M$.
\par
Thus we have a nice description of $Simp_{2}(S)$ as the
open subset of the affine 5-space $Spec(C_{2})$, defined by $d_{[X,Y]}\neq 
0$. It should be compared with the computations in section 
\ref{planesimp2}.

\subsection{Simple 2-dimensional modules for noncommutative plane
curves}

In classical algebraic geometry the simple finite dimensional
representations of a commutative $k$-algebra $A_{0}=k[x,y]/(f_{0})$ are 
parametrised by a
plane curve given by the zero set of the polynomial $f_{0}$.
A noncommutative model of the plane curve $f$,
is a (noncommutative) $k$-algebra $A=k\langle x,y\rangle /(f)$
where $f\equiv f_{0}\,([x,y])$. The
1-dimensional representations of $A$ are the same as for the commutative
curve $A_{0}$. But, in contrast to what is the case for a 
commutative algebra, $A$ may also have higher 
dimensional simple representations. 
\par
The purpose of this section is to give
a description of the 2-dimensional simple representations of 
noncommutative plane curves.
\begin{Lem}\label{traceformulas}
    Let $R$ be a commutative $k$-algebra.
    For any two $2\times 2$-matrices $X,Y\in M_{2}(R)$ we have the 
    equalities
    \begin{itemize}
	\item[i)] $X^2=t_{X}X-d_{X}$
	\item[ii)] $X^3=(t_{X}^2-d_{X})X-t_{X}d_{X}$
	\item[iii)] $YX=-XY+t_{X}Y+t_{Y}X+t_{XY}-t_{X}t_{Y}$
    \end{itemize}
\end{Lem}
\begin{proof}
    Equality i) and ii) follows from the Cayley-Hamilton theorem. 
    To prove iii), notice that
    \begin{displaymath}
	(X+Y)^2=t_{X+Y}(X+Y)-d_{X+Y}I
    \end{displaymath}
    On the other hand we have
    \begin{displaymath}
	(X+Y)^2=X^2+XY+YX+Y^2=t_{X}X-d_{X}+XY+YX+t_{Y}Y-d_{Y}
    \end{displaymath}
    Thus, using the fact that $d_{X+Y}-d_{X}-d_{Y}=t_{X}t_{Y}-t_{XY}$
    iii) follows.
\end{proof}
In section \ref{rep} we defined a versal family
\begin{displaymath}
    \Phi:k\langle x,y\rangle \rightarrow M_{2}(\Gamma)
\end{displaymath}
where $\Gamma=k[x_{p,q}^i]/(\tilde{f})$.
Let
$X:=\Phi(x)$ and $Y:=\Phi(y)$. Then for $f\in k\langle x,y\rangle$ we 
have
$\Phi(f)=c_{1}XY+c_{2}X+c_{3}Y+c_{4}I$, where the coefficients 
$c_{i}\in \Gamma$, $i=1,\ldots,4$ are uniquely determined by $f$.
This follows from the Cayley-Hamilton theorem and from the linear 
independence of the set $\{I,X,Y,XY\}$. 
\par
The simplicity criterion for a representation is given by
the non-vanishing of the Formanek element
\begin{displaymath}
    d_{[X,Y]}=
    -\frac{1}{4}((2t_{XY}-t_{X}t_{Y})^2-(t_{X}^2-4d_{X})(t_{Y}^2-4d_{Y}))
\end{displaymath}
Now, let $\phi:k\langle x,y\rangle\rightarrow M_{2}(k)$
be a 2-dimensional representation of 
$k\langle x,y\rangle$, corresponding to a closed point $\gamma\in 
Spec(\Gamma)$. Then $\phi$ induces a simple
representation of $A=k\langle x,y\rangle/(f)$ if and only if
$c_{i}(\gamma)=0$ for $i=1,\ldots,4$ and $d_{[X,Y]}(\gamma)\not=0$.
\par
\textbf{Example}. For $\delta\in k$ we put $f_{\delta}=x^2+y^2-1+\delta [x,y]$ and 
$A_{\delta}=k\langle x,y\rangle/(f_{\delta})$. Using lemma 
\ref{traceformulas} we find
\begin{align*}
    \Phi(x^2+y^2&-1+\delta[x,y])= \\
    2\delta XY + &(t_{X}-\delta t_{Y})X
    +(t_{Y}-\delta t_{X}) Y + (-d_{X}-d_{Y}-1-\delta t_{XY}+\delta
    t_{X}t_{Y})\\
\end{align*}
Thus $Simp_{2}(A_{\delta})=\emptyset$ 
unless $\delta=0$. For $\delta=0$ we have
\begin{displaymath}
    Spec(\Gamma(A_{\delta=0}))
    =V(t_{X},t_{Y},d_{X}+d_{Y}+1)\subset {\mathbf A}^4
\end{displaymath}
and 
$Simp_{2}(A_{\delta=0})$ is the open subscheme of 
\begin{displaymath}
    Spec(\Gamma(A_{\delta=0}))\subset {\mathbf A}^4
\end{displaymath}
given by the non-vanishing of the Formanek element 
\begin{displaymath}
      t_{XY}^2-4d_{X}(-1-d_{X})=t_{XY}^2+(2d_{X}+1)^2-1
\end{displaymath}

\subsection{The general case}\label{generalcase}

Now for a general finitely generated algebra $A=S/I,\ I=(f_1,..,f_r)$
we may, as
above, consider the trace ring $C_n\subseteq \Gamma^G$. One might hope that
$Spec(C_n)=Spec(\Gamma^{G})$, and that this space parametrises the
semisimple $n$-dimensional representations. However, this is in 
general not true,
as we shall see in example \ref{exampledualnumbers}.
\par
In the paper (\cite{DP}) De Concini and Procesi handle 
the problem by
restricting to algebras with trace. This limitation
excludes some interesting examples and weakens the results 
substancially.
\par
One of the weak points of this approach is related to the embedding of
$Simp_{n} (A)$ in the
representation space. The closure of $Simp_{n} (A)$ in this space consist of
semi-simple modules, which ``do not deform back into simple modules''.
\par
We feel that
this picture should be extended to include completions of
$Simp_n(A)$ and therefore include the indecomposable modules, see section
\ref{completions}.
To achieve this we need to take a different point of view. 

\subsection{Example}\label{exampledualnumbers}

The set of 2-dimensional representations of the polynomial ring 
$S=k[x]$ is parametrised by the variety $\X_{S}^{(2)}=Spec(\Gamma_{S})$ where 
\begin{displaymath}
    \Gamma_{S}=k[x_{11},x_{12},x_{21},x_{22}]
\end{displaymath}
\par
Let $A=k[x]/(x^2)$. Then
\begin{displaymath}
    \X=\X_{A}^{(2)}=Spec(\Gamma)
\end{displaymath}
where, as we have seen in section \ref{rep}, 
\begin{displaymath}
    \Gamma=k[x_{11},x_{12},x_{21},x_{22}]/J
\end{displaymath}
and $J$ is the two-sided ideal
\begin{displaymath}
    J=(x_{11}^2+x_{12}x_{21},(x_{11}+x_{22})x_{12},
    (x_{11}+x_{22})x_{21},x_{22}^2+x_{12}x_{21})
\end{displaymath}
Thus $\X$ is a 
2-dimensional subvariety of $\mathbf{A}^4=Spec (\Gamma_{S})$. The 
versal family is given by the subring of the image of the map
\begin{displaymath}
    \pi:A\longrightarrow \Gamma\otimes M_{2}(k)
\end{displaymath}
given by $x\mapsto \sum_{i,j}\overline{x}_{ij}\otimes e_{ij}$. The 
trace ring $C_{2}$ is generated by the traces of the matrices in the 
versal family. It is easily seen that $C_{2}\simeq k[t]$ where 
$t=\overline{x}_{11}+\overline{x}_{22}$. We call 
$Spec(C_{2})=\mathbf{A}^1$ the trace space. 
\par
The action of the group $G=PGl_{2}$ on $\X$ has two orbits, the 
trivial semi-simple module $V_{0}$, corresponding to the origin, and 
$V_{1}$, corresponding to the indecomposable module given by 
$x\mapsto \2matr 0100$. Thus $\X/G$ consists of two points, one in 
the closure of the other.
\par
The ideal $J$ is $G$-invariant and $G$ acts on the short-exact sequence 
\begin{displaymath}
    0\rightarrow J\longrightarrow \Gamma_{S}\longrightarrow 
    \Gamma\rightarrow 0
\end{displaymath}
inducing an long-exact sequence 
\begin{displaymath}
    0\rightarrow H^0(G,J)\longrightarrow H^0(G,\Gamma_{S})
    \longrightarrow H^0(G,\Gamma)\longrightarrow H^1(G,J)\rightarrow \ldots
\end{displaymath}
But $H^1(G,J)=0$ since $J$ is a graded $k$-vector space and the 
reductive group $G$ acts 
on each finite dimensional graded component. Therefore the sequence 
reduces to a short-exact sequence
\begin{displaymath}
    0\rightarrow J\cap \Gamma_{S}^{G}\longrightarrow \Gamma_{S}^{G}
    \longrightarrow \Gamma^{G}\rightarrow 0
\end{displaymath}
where $\Gamma_{S}^{G}=k[x_{11}+x_{22},x_{11}x_{22}-x_{12}x_{21}]$.
Since $J\cap \Gamma_{S}^{G}$ is precisely the ideal of 
$\Gamma_{S}^{G}$ generated by the determinant of the generic matrix, 
we have
\begin{displaymath}
    \Gamma^{G}\simeq k[t]
\end{displaymath}
and consequently $Spec(\Gamma^{G})\not\simeq Spec(\Gamma)/G$.
\par
In section \ref{planesimp2} we shall discuss this example a bit further.

\subsection{Local approach using deformation theory}\label{local}

A different approach to understanding the geometry of $Simp_{n}(A)$ 
is to analyse the local structure of $Simp_n(A)$
using deformation theory, see \cite{La3.1},\cite{La1}
of which we briefely recall the main content. The formal moduli of an
$n$-dimensional simple $A$-module $V$ is a complete (noncommutative)
$k$-algebra $H^A(V)$ (\cite{La3.1}). Using the structure map of $V$
and the basic properties of the formal moduli one can show that there
is a natural,
topological surjective
$k$-algebra homomorphism (see \cite{La4})
\begin{displaymath}
    A\rightarrow H^A(V)\otimes_{k} \textrm{End}_{k}(V)\simeq 
    M_{n}(H^A(V))
\end{displaymath}
Recall that a standard $n$-commutator relation in a $k$-algebra $A$ is
a relation of the type,
\begin{displaymath}
    [a_1, a_2, ..., a_{2n}]:=\sum_{\sigma\in \Sigma_{2n}} sign(\sigma)
    a_{\sigma(1)} a_{\sigma(2)} ... a_{\sigma(2n)} =0
\end{displaymath}
\noindent where $\{a_1, a_2, ..., a_{2n}\}$ is a subset of $A$. Let $I(n)$ be
the two-sided ideal of $A$ generated by the subset,
\begin{displaymath}
    \{[a_1, a_2, ..., a_{2n}]|\ \{a_1, a_2, ..., a_{2n}\}\subset A\}.
\end{displaymath}
\noindent Any $k$-algebra homomorphism $\phi:A\rightarrow M_{n}(k)$ factors
through $A(n):=A/I(n)$, and $Simp_{m}(A)= Simp_{m}(A(n))$ for $m\le
n$, see e.g. \cite {F5}. For $m>n$ we have $Simp_{m}(A(n))=\emptyset$.
Notice that $A(1)$ is the commutativisation of $A$.\\
\par
A usefull fact is that the formal
moduli of $V\in Simp_{n}(A)$ considered as an $A(n)$-module is 
isomorphic to the commutativisation of
the formal moduli of $V$ over $A$;
\begin{displaymath}
    H^{A(n)}(V)\simeq H^A(V)_{\textrm{com}}
\end{displaymath}
\noindent Consequently $H^{A(n)}(V)$ is a commutative
$k$-algebra, easily computed when we know $H^{A}(V)$.
\par
For a free algebra $S=k\langle x_{1},\dots,x_{m}\rangle$ on $m$
variables and $V\in Simp_{n}(S)$ the formal
moduli is a noncommutative power series algebra generated by
$m\cdot n^2-(n^2-1)$ elements, thus in that case we have
\begin{displaymath}
    H^{S(n)}(V)\simeq k[[t_{1}, \dots,t_{(m-1)n^2+1}]]
\end{displaymath}
\noindent which should be compared to the Theorem \ref{proc}
of Procesi given above.
\par
Now consider the product morphism
\begin{displaymath}
    \mu: A(n)\rightarrow \prod_{V\in Simp_n(A)}H^{A(n)}(V)\otimes_k End_k(V)
\end{displaymath}
\noindent In general this map
is not injective. Put 
\begin{displaymath}
    B=\prod_{V\in Simp_n(A)}H^{A(n)}(V)
\end{displaymath}
Let $x_i\in A$, $i=1,\ldots,m$ be generators of $A$,
and consider the images 
\begin{displaymath}
    \mu(x_{i})=(x^i_{p,q})\in B\otimes_k End_k(k^n)\simeq M_{n}(B)
\end{displaymath}
obtained by choosing bases in all $V\in Simp_n(A)$. Now, 
since $B$ is
commutative, the $k$-subalgebra $C(n)\subset B$
generated by the elements $\{x_{p,q}^i\}_{i=1,..,d;\ p,q=1,..,n}$ is
commutative.
One can show (see \cite{La4}) that there is an open embedding 
$Simp_n(A)\subset Simp_{n}(C(n))$.
Thus any simple, $n$-dimensional $A$-module corresponds
to a closed point $v\in Simp_{1}(C(n))$.
\par
The space $Simp_{1}(C(n))$ is in a rather straightforward sense, a {\it
compactification} of $Simp_n(A)$, as is $Simp_{1}(C_{n})$, but neither 
are good completions of $Simp_{n}(A)$.
\par
The complement of the simple
locus of $Simp_{1}(C_{n})$ is in \cite{Pr2} identified with the set of 
semisimple modules. 
However, for a non-simple, semisimple $A$-module $M$
we have 
\begin{displaymath}
    dim_{k}(End_{A}(M))\ge 2
\end{displaymath}
while for an indecomposable module $E$
which is an extension of non-isomorphic 1-dimensional representations 
$End_{A}(E)=k$. Any modular
infinitesimal deformation preserves the dimension of the endomorphism 
ring, see \cite{La3.1}.
Therefore there is no modular deformation of a decomposable module $M$ 
into a simple module, i.e. a 
non-simple semi-simple module $M$ is not in the completion of any 
$Simp_{n}(A)$.
\par
Let us return to our example from the previous section.
\par
\textbf{Continuation of example \ref{exampledualnumbers}}. Recall that we
considered the 2-dimensional representations 
of the $k$-algebra $A=k[x]/(x^2)$. The tangent space of the modular 
part of the deformations of the semi-simple representation $V_{0}$ is given by 
\begin{align*}
    Ext_{A}^1(V_{0},V_{0})^{End_{A}(V_{0})}&=
    \{\xi\in Ext_{A}^1(V_{0},V_{0})\,\vert\,
    \phi^{\ast}\xi=\xi\phi_{\ast}\,\,
    \forall \phi\in End_{A}(V_{0})\}\simeq k
\end{align*}
Since the cup product $\xi\cup\xi\in Ext_{A}^2(V_{0},V_{0})$ 
is non-zero the formal moduli is given by $k[y]/(y^2)$, where $y$ corresponds 
to the tangent direction $\xi$. So the moduli space $\X/G$ consists 
of two points, the semi-simple module $V_{0}$, and  
the indecomposable $V_{1}$. The $G$-orbit of $V_{1}$ has $V_{0}$ in its closure.

\subsection{The free noncommutative algebra in two variables, again}
\label{planesimp2}

Let us once more consider the noncommutative affine plane
$S=k\langle x,y\rangle $ and its 2-dimensional representations. This 
time we use the local approach via deformation theory.
\par
In \cite{La4} Laudal has computed a versal family of 
2-dimensioanl representations of $S$, as deformations of the simple 
module $V\in Simp_{2}(S)$, given by 
\begin{displaymath}
    x\mapsto \2matr 0100, \qquad y\mapsto \2matr 0010
\end{displaymath}
The formal moduli is the smooth, formal polynomial 
$k$-algebra 
\begin{displaymath}
    H^{S}(V)_{com}\simeq k[[t_{1},t_{2},t_{3},t_{4},t_{5}]]
\end{displaymath}
and the versal family is defined by 
\begin{displaymath}
    x\mapsto \2matr 0{1+t_{3}}{t_{5}}{t_{4}}, \qquad 
    y\mapsto \2matr {t_{1}}{t_{2}}{1+t_{3}}0
\end{displaymath}
The algebra $C(2)$ is therefore given by 
\begin{displaymath}
    C(2)\simeq  k[t_{1},t_{2},t_{3},t_{4},t_{5}]
\end{displaymath}
Clearly $C_{2}\varsubsetneq C(2)$, where $C_{2}$, given in section 
\ref{planesimp}, is generated by
\begin{align*}
    t_{X}&=t_{4}, \qquad t_{Y}=t_{1}\\
    d_{X}&=-t_{5}(1+t_{3}), \qquad d_{Y}=-t_{2}(1+t_{3})\\
    t_{XY}&=(1+t_{3})^2+t_{2}t_{5}
\end{align*}
with Formanek element
\begin{displaymath}
    d_{[X,Y]}=-((1+t_{3})^2-t_{2}t_{5})^2
    +(t_{1}(1+t_{3})+t_{2}t_{4})(t_{4}(1+t_{3})+t_{1}t_{5})
\end{displaymath}
Notice that the indecomposable, non-simple modules obtained by e.g. 
letting $t_{5}=1+t_{3}=0$, are not visible in $Spec(C_{2})$. 
The family of non-isomorphic indecomposable modules, parametrised by 
$t_{2}$, collapses in the trace ring to one semi-simple module.


\section{Extensions}\label{completions}

The simple $A$-modules of dimension $n$ are parametrised by an open set
in some affine
space. One of the main objectives of this paper is to construct 
natural completions
of the simple locus. The completions are formed by adding indecomposable modules 
obtained as iterated extensions of simple modules
of lower dimension.

\subsection{The Noncommutative Jacobi matrix}\label{Jacobi}

Let $S=k\langle x_{1,}\ldots,x_{m}\rangle$ be a free $k$-algebra on
$m$ non-commuting variables. Let $f\in S$ be a polynomial. Let
$\phi_{p}:S\rightarrow k(p)$ be a 1-dimensional representation 
corresponding to a point $p\in Simp_{1}(S)\simeq \mathbf{A}^m$, such that
$\phi_{p}(f)=f(p)=0$. The representation is given by $\phi_{p}(x_{i})=a_{i}$ 
for some
$a_{i}\in k$ and $i=1,2,\ldots,m$, i.e. $p=(a_{1},\ldots,a_{m})\in 
\mathbf{A}^{m}$.
\begin{Defn}
    Let $f\in S$ and let $\phi_{p}:S\rightarrow k(p)$ be a 
    1-dimensional representation, corresponding to the point 
    $p=(a_{1},\ldots,a_{m})\in \mathbf{A}^{m}$.
    An equality
    \begin{displaymath}
	f=\sum_{k=1}^m f_{k,p}(x_{k}-a_{k})
    \end{displaymath}
    is called a left decomposition of $f$ with respect to $k(p)$.
\end{Defn}
\begin{Prop}
    For any $f\in S$ and any $\phi_{p}:S\rightarrow k(p)$ satisfying $\phi(f)=0$
    there exists a unique left decomposition of $f$ with respect to
    $k(p)$.
\end{Prop}
\begin{proof}
    For any monomial $\underline{m}$ in $f$ of positiv degree, we can write
    \begin{displaymath}
	\underline{m}=\underline{m}^{\prime}(x_{i}-a_{i})
	+a_{i}\underline{m}^{\prime}
    \end{displaymath}
    for suitable $a_{i}\in k$, and where $\underline{m}^{\prime}$ is
    of strictly less degree than $\underline{m}$. 
    Thus by an inductive procedure we can write
    \begin{displaymath}
	f(\underline{x})=f(\underline{a})+
	\sum_{k=1}^m f_{k,p}(x_{k}-a_{k})
    \end{displaymath}
    The existence of a left decomposition follows from the fact that
    $f(\underline{a})=\phi_{p}(f)=0$.
    \par
    Since the set $\{x_{k}-a_{k}\}_{k=1,\ldots,m}$ is
    a free generating set for $k\langle x_{1},\ldots,x_{m}\rangle$ 
    the decomposition is obviously unique.
\end{proof}
There exists an algorithm for computing the components of this decomposition.
Let $\phi_{p}:S\rightarrow k(p)$ and let $D_{i}(-;p)$ be the linear form 
defined on $S$ such that $D_{i}(a;p)=0$ for $a\in k$,
$D_{i}(x_{j};p)=\delta_{ij}$ and for a product $fg$
\begin{displaymath}
    D_{i}(fg;p)=fD_{i}(g;p)+D_{i}(f;p)g(p)
\end{displaymath}
So $D_{i}(-;p)$ act as a partial derivation, satisfying the Leibniz
rule where we evaluate in $p$ on the right hand side. In fact 
\begin{displaymath}
    D_{i}(-;p)\in Der_{k}(S,Hom_{k}(k(p),S)),
\end{displaymath}
which proves that $D_{i}(\textrm{-};p)$ is well-defined.
\begin{Defn}
    The element $D_{k}(f;p)$ is called the noncommutative
    (left) $k$-th partial derivative
    of $f$ with respect to the 1-dimensional representation $k(p)$.
\end{Defn}
\begin{Prop}\label{firstorder}
    The left components of an element $f\in S$ with respect to 
    $k(p)$, where $p=(a_{1},\ldots,a_{m})$, is precisely the left
    partial derivatives of $f$, i.e. the equality
    \begin{displaymath}
	f=\sum_{k=1}^m D_{k}(f;p)(x_{k}-a_{k})
    \end{displaymath}
    holds.
\end{Prop}
\begin{proof}
    If one applies $D_{k}(-;p)$ to the equality
    \begin{displaymath}
	f=\sum_{k=1}^m f_{k,p}(x_{k}-a_{k})
    \end{displaymath}
    and use the Leibniz rule, then it follows that $f_{k,p}=D_{k}(f;p)$.
\end{proof}
We can extend the definition of noncommutative partial derivative
to a two-sided ideal $I\subset S$. Let $I\subset S$ be a two-sided 
ideal generated by 
$f^{1},\ldots,f^{r}$ and $\phi_{p}$ be a
representation of $S$ corresponding to a point $p\in Simp_{1}(S/I)$. 
Consider the left ideal of $A=S/I$
generated by the image of the $i$-th partial derivatives of the
generators $f^{1},\ldots,f^{r}$. Denote this ideal
$J_{i}(I;f^{1},\ldots,f^{r};p)$.
\begin{Lem}
    The ideal $J_{i}(I;f^{1},\ldots,f^{r};p)$ is independent of the
    choice of generator set for the ideal $I$.
\end{Lem}
\begin{proof}
    Let
    \begin{displaymath}
	f=\sum_{j,k}\alpha_{jk}f^j\beta_{jk}
    \end{displaymath}
    be an arbitrary element of $I$, with $\alpha_{jk},\beta_{jk}\in S$.
    Then we have
    \begin{align*}
	D_{i}(f;p)
	&=D_{i}(\sum_{j,k}\alpha_{jk}f^j\beta_{jk};p)\\
	&=\sum_{j,k}(D_{i}(\alpha_{jk};p)\phi_p(f^j\beta_{jk})
	+\alpha_{jk}D_{i}(f^j;p)\phi_p(\beta_{jk})
	+\alpha_{jk}f^jD_{i}(\beta_{jk};p)\\
	&\in J_{i}(I;f^{1},\ldots,f^{m};p) + I
    \end{align*}
\end{proof}
We write $J_{i}(I;p)$ for this ideal.
\begin{Defn}
    The ideal $J_{k}(I;p)\subset S$ is called the
    $k$-th noncommutative Jacobi ideal of $I$ with respect to $k(p)$.
    The matrix $J(I;p)=(J_{k}(f^j;p))_{k,j}$ is called the
    noncommutative Jacobi matrix
    of the presentation $I=(f^1,\ldots,f^r)$.
\end{Defn}
Notice that for any element $g=\sum_{i}r_{i}(x_{i}-a_{i})\in S$ we have
$r_{i}=D_{i}(g;p)$ by the uniqueness of partial derivatives. If $g\in
I$ the proof of the lemma shows that
\begin{displaymath}
    r_{i}=\sum_{j}(\sum_{k}\alpha_{jk}\phi_p(\beta_{jk}))D_{i}(f^j;p)
    +\sum_{j,k}\alpha_{jk}f^jD_{k}(\beta_{jk};p)
\end{displaymath}
Modulo the ideal $I$ we can write
\begin{displaymath}
    r_{i}=\sum_{j}s_{j}D_{i}(f^j;p)
\end{displaymath}
where
\begin{displaymath}
    s_{j}=\sum_{k}\alpha_{jk}\phi_p(\beta_{jk})
\end{displaymath}
A consequence of this is that
\begin{displaymath}
    A^r\buildrel{J(I;p)}\over{\longrightarrow}A^m
    \buildrel{(\underline{x}-\underline{a})}\over{\longrightarrow}
    A\buildrel{\phi_p}\over{\longrightarrow}k(p)\rightarrow 0
\end{displaymath}
gives a truncated resolution of the representation $k(p)$ as a left
$A$-module.
\par
The noncommutative Jacobi matrix $J(I;p)$ is an $r\times m$-matrix 
over $S$. For a representation $\phi_{2}:S\rightarrow k(p_{2})$ we 
write $J(I;p_{1})(p_{2})$ for the evaluation of $J(I;p_{1})$ in 
$p_{2}$.
\begin{Prop}\label{extdim}
    Let $A=k\langle x_{1,}\ldots,x_{m}\rangle/I$ be a $k$-algebra and let 
    $\phi_{1},\phi_{2}$ be 1-dimensional representations corresponding to 
    $k(p_{1}),k(p_{2})\in {\Simp}_{1}(A)$. Suppose $p_{1}\not=p_{2}$, 
    then
    \begin{displaymath}
	\textrm{dim}_{k} \ext{1}{A}{k(p_{1})}{k(p_{2})}
	=m-1-\textrm{rk}\,J(I;p_{1})(p_{2}).
    \end{displaymath}
    If $p_{1}=p_{2}=:p$, then 
    \begin{displaymath}
	\textrm{dim}_{k} \ext{1}{A}{k(p)}{k(p)}
	=m-\textrm{rk}\,J(I;p)(p).
    \end{displaymath}
    
\end{Prop}
\begin{proof}
    A direct consequence of the above discussion.
\end{proof}   
For $A=S/I$ commutative the commutators $x_{i}x_{j}-x_{j}x_{i}$ are part
of a generator set for the ideal $I\subset S$. Notice that for 
$p=(a_{1},\ldots,a_{m})$
\begin{displaymath}
    D_{k}(x_{i}x_{j}-x_{j}x_{i};p)=
    \begin{cases}
	0&k\not=i,j\\
	a_{j}-x_{j}&k=i\\
	x_{i}-a_{i}&k=j\\
    \end{cases}
\end{displaymath}
Thus $J_{k}(I;p_{1})\subset \textbf{m}$, the maximal ideal of $S$
corresponding to $p_{1}$. Therefore, for $p_{1}\not=p_{2}$ we have
\begin{displaymath}
    \textrm{rk}\,J_{k}(I;p_{1})(p_{2})=m-1
\end{displaymath}
proving the well-known vanishing of $\ext{1}{A}{k(p_{1})}{k(p_{2})}$ for $A$
commutative and $p_{1}\not=p_{2}$. For $A$ noncommutative this
is no longer  true, and a natural question to ask is, given a 
representation $\phi_p$ of $A$, corresponding to a point 
$k(p)\in{\Simp}_{1}(A)$, describe the variety of simple 1-dimensional
representations $\phi_{q}$ with non-vanishing $Ext^1_{A}(k(p),k(q))$. Under
certain finiteness conditions this gives a correspondance on the set
of simple 1-dimensional representations. In the next section we shall
explore this in the case of plane curves.

\subsection{Noncommutative Taylor series expansion}\label{Taylor}

What we did in the previous section is a special case of a more 
general set-up, see \cite{La4}. Let $m$ be some natural number
and consider the free $k$-algebra on $3m$ variables,
\begin{displaymath}
    k\langle \underline{x},\underline{v},\underline{u}\rangle=
    k\langle x_{1},\ldots,x_{m}, v_{1},\ldots,v_{m},u_{1},\ldots,u_{m}
    \rangle
\end{displaymath}
\begin{Defn}
    Let $f\in S=k\langle x_{1},\ldots,x_{m}\rangle$. Denote by
    \begin{displaymath}
	D_{x_{i}}(f;\underline{x},\underline{v},\underline{u})\in 
	k\langle \underline{x},\underline{v},\underline{u}\rangle
    \end{displaymath}
    the linear function in $f$, defined for $f=x_{i}$ and 
    $f=\underline{m}x_{j}$, $\underline{m}$ any monomial,  
    by
    \begin{align*}
	D_{x_{i}}(x_{j};\underline{x},\underline{v},\underline{u})
	&=\delta_{i,j}v_{i}\\
	D_{x_{i}}(\underline{m}x_{j};\underline{x},\underline{v},\underline{u})
	&=D_{x_{i}}(\underline{m};\underline{x},\underline{v},\underline{u})u_{j}
	+\delta_{i,j}\underline{m}v_{i}\\
    \end{align*}
\end{Defn}
The $D_{x_{i}}(f)$ is a general noncommutative differentiation 
symbol. We can as well define higher order differentiation by 
iteration of the function $D$, considering the 
$(\underline{v},\underline{u})$-variables as constants.
Inductively we use the notation 
\begin{displaymath}
    D_{x_{i_{n}},\ldots,x_{i_{2}},x_{i_{1}}}
    (f;\underline{x},\underline{v},\underline{u})
    =D_{x_{i_{n}}}(D_{x_{i_{n-1}},\ldots,x_{i_{1}}}
    (f;\underline{x},\underline{v},\underline{u});
    \underline{x},\underline{v},\underline{u})
\end{displaymath}
\begin{Prop}
    For $f\in S$ we have the equality in $k\langle 
    \underline{u},\underline{v}\rangle$, given by
    \begin{displaymath}
	f(\underline{u}+\underline{v})=f(\underline{u})
	+\sum_{n}\sum_{1\le i_{1},\ldots,i_{n}\le d}
	D_{x_{i_{n}},\ldots,x_{i_{2}},x_{i_{1}}} 
	(f;\underline{u},\underline{v},\underline{u})
    \end{displaymath}
\end{Prop}
\begin{proof}
    Reduce to the case where $f$ is a monomial and use a 
    straightforward combinatorial argument, see \cite{La4}.
\end{proof}
If we put $\underline{v}=\underline{x}-\underline{u}$ into this 
formula we get a noncommutative Taylor expansion
\begin{displaymath}
    f(\underline{x})=f(\underline{u})
    +\sum_{n}\sum_{(i_{1},\ldots,i_{n})\in\{1,\ldots,m\}^n}
    D_{x_{i_{n}},\ldots,x_{i_{2}},x_{i_{1}}} 
    (f;\underline{u},\underline{x}-\underline{u},\underline{u})
\end{displaymath}   
Notice that if the variable set $\underline{u}$ is assumed to be in 
the center of the algebra, i.e. if we work over 
$k[\underline{u}]\langle\underline{x}\rangle$, the 
Taylor expansion takes the more familiar form,
\begin{displaymath}
    f(\underline{x})=f(\underline{u})
    +\sum_{n}\sum_{1\le i_{1},\ldots,i_{n}\le d}
    D_{x_{i_{n}},\ldots,x_{i_{2}},x_{i_{1}}} 
    (f;\underline{u},\underline{u})
    (x_{i_{n}}-u_{i_{n}})\ldots(x_{i_{1}}-u_{i_{1}})
\end{displaymath}
where we use the notation 
\begin{displaymath}
    D_{x_{i_{n}},\ldots,x_{i_{2}},x_{i_{1}}} 
    (f;\underline{x},\underline{u})
    =D_{x_{i_{n}},\ldots,x_{i_{2}},x_{i_{1}}} 
    (f;\underline{x},(1,\ldots,1),\underline{u})
\end{displaymath}
The connection with $D_{k}(f;p)$ as defined in the previous 
section is given by the next proposition.
\begin{Prop}\label{firstorder2}
    For any $f\in k\langle \underline{x}\rangle$ we have the following 
    equality in the extended $k$-algebra 
    $k[\underline{u}]\langle\underline{x}\rangle$,
    \begin{displaymath}
	f(\underline{x})
	=f(\underline{u})
	+\sum_{i=1}^m D_{x_{i}}(f;\underline{x},\underline{u})(x_{i}-u_{i})
    \end{displaymath}
\end{Prop}
\begin{proof}
    Applying the Taylor formula to the polynomial 
    $D_{x_{i_{1}}}(f;\underline{x},\underline{u})$ we get the equality
    \begin{align*}
	D_{x_{i_{1}}}&(f;\underline{x},\underline{u})
	=D_{x_{i_{1}}}(f;\underline{u},\underline{u})\\
	+&\sum_{n}\sum_{1\le i_{2},\ldots,i_{n}\le d}
	D_{x_{i_{n}},\ldots,x_{i_{2}},x_{i_{2}}} 
	(D_{x_{i_{1}}}(f;\underline{u},\underline{u});\underline{u},\underline{u})
	(x_{i_{n}}-u_{i_{n}})\ldots(x_{i_{2}}-u_{i_{2}})
    \end{align*}
    Multiplying from right by $(x_{i_{1}}-u_{i_{1}})$, adding up over 
    all $i=1,2,\ldots,d$ and using the Taylor formula once 
    more gives the result.
\end{proof}
Since the left decomposition is unique we it follows that 
\begin{displaymath}
    D_{k}(f;p)=D_{x_{k}}(f;\underline{x},p)
\end{displaymath}

\subsection{Extension relations on plane curves}\label{alg.cor}

Let $A=k \langle x,y\rangle /(f)$ be a noncommutative model for 
the reduced algebraic plane curve,
\begin{displaymath}
    C=Spec(k[x,y]/(f_{0})):=Simp_1(A)
\end{displaymath}
and define a relation $\R$ on the affine plane $\mathbf{A}^2$ by
\begin{displaymath}
    \R=V(D_{1}(f;\underline{u})(\underline{v}),
    D_{2}(f;\underline{u})(\underline{v}))\subset \mathbf{A}^2\times 
    \mathbf{A}^2
\end{displaymath}
where the product is parametrised by $(\underline{u},\underline{v})$.
\begin{Thm}
    Let $A=k\langle x,y\rangle/(f)$ and $C=Simp_{1}(A)=Simp_{1}(A_{0})$.
    \begin{itemize}
	\item[i)] The restriction of $\R$ to $C$ induces a relation on $C$, 
	i.e.
	\begin{displaymath}
	    \R\cap(C\times\mathbf{A}^2)=C\times C
	\end{displaymath}
	Moreover, $\R\cap\Delta C=Sing(C)$, the singular locus of $C$.
	\item[ii)] Let $p_{1},p_{2}\in C$ be two points 
	corresponding to the simple modules $k(p_{1})$ and 
	$k(p_{2})$. Then 
	there exist a noncommutative model 
	$A=k\langle x,y\rangle/(\tilde{f})$ of $f_{0}$ such that 
	$Ext_{A}^1(k(p_{1}),k(p_{2}))\neq 0$.
	\par
	\item[iii)] Consider models given by  $\tilde{f}_{t}=f+t[x,y]$. Then for 
	generic $t$ the projection $pr_{1}:\R\cap(C\times C)\rightarrow C$ is 
	dominant.
    \end{itemize}
\end{Thm}
\begin{proof}
    i) By Prop. \ref{firstorder} we have for $p=(\alpha,\beta)\in C$
    \begin{displaymath}
	f=D_{1}(f;p)(x-\alpha)+D_{2}(f;p)(y-\beta).
    \end{displaymath}
    If for $q\in\mathbf{A}^2$ we have 
    \begin{displaymath}
	D_{1}(f;p)(q)=D_{2}(f;p)(q)=0
    \end{displaymath}
    Then this implies $f(q)=0$, i.e. $q\in C$.
    \par
    Now use Prop. \ref{extdim} and see that for $p=q$ the above 
    equations imply $J((f);p)(p)$ has rank zero. So $p\in C$ is a 
    singular point.
    \par
    ii) Let $\tilde{f}=f+t[x,y]$ for some $t\in k$ and 
    let $p_{1},p_{2}\in C$ be two points corresponding to the simple 
    modules $k(p_{1})$ and $k(p_{2})$. By Prop. \ref{extdim} 
    $Ext_{A}^1(k(p_{1}),k(p_{2}))\neq 0$ if and only if there 
    exists a scalar $t$ such that 
    \begin{align*}
	D_{1}(\tilde{f};p_{1})(p_{2})
	&=D_{1}(f;p_{1})(p_{2})+t\cdot (b_{1}-b_{2})=0\\
	D_{2}(\tilde{f};p_{1})(p_{2})
	&=D_{2}(f;p_{1})(p_{2})+t\cdot (a_{2}-a_{1})=0\\
    \end{align*}
    or equivalently 
    \begin{displaymath}
	D_{1}(f;p_{1})(p_{2})(a_{2}-a_{1})+D_{2}(f;p_{1})(p_{2})(b_{2}-b_{1})=0
    \end{displaymath}
    Using Prop. \ref{firstorder} we identify the left hand side of 
    this equality as $f(p_{2})-f(p_{1})$ and the result follows.
    \par
    iii) Since $\R\cap\Delta =Sing(C)$ is finite, we just need to 
    show that the $k$-algebra homomorphism
    \begin{displaymath}
	pr_{1}^{\ast}:k[u_1,u_2]/(f_{0})\longrightarrow 
	k[u_1,u_2,v_{1},v_{2}]/
	(D_{1}(\tilde{f}_{t};\underline{u})(\underline{v}),
	D_{2}(\tilde{f}_{t};\underline{u})(\underline{v}),f_{0}(\underline{u})),
	f_{0}(\underline{v}))
    \end{displaymath}
    is injective. Suppose it is not injective. Then there must exist 
    \begin{displaymath}
	h_{1},h_{2}\in k[u_1,u_2,v_{1},v_{2}]
    \end{displaymath}
    such that 
    \begin{align*}
	h_{1}D_{1}(\tilde{f}_{t};\underline{u})(\underline{v})&+
	h_{2}D_{2}(\tilde{f}_{t};\underline{u})(\underline{v})\\
	&=h_{1}(D_{1}(f;\underline{u})(\underline{v})+t(u_{2}-v_{2}))+
	h_{2}(D_{2}(f;\underline{u})(\underline{v})+t(v_{1}-u_{1}))\\
	&=h_{1}D_{1}(f;\underline{u})(\underline{v})
	+h_{2}D_{2}(f;\underline{u})(\underline{v})
	+t(h_{1}(u_{2}-v_{2})+h_{2}(v_{1}-u_{1}))\\
	&\in (u_1,u_2)
    \end{align*}
    For this to be true for generic $t$, we must have  
    \begin{displaymath}
	h_{1}(u_{2}-v_{2})+h_{2}(v_{1}-u_{1})\in (u_1,u_2)
    \end{displaymath}
   The indeterminants $v_{1}$ and $v_{2}$ are linearly independent over 
   $k[u_1,u_2]$, hence there must exist $h\in 
   k[\underline{u},\underline{v}]$ such thaht  
   $h_{1}=h\cdot (v_{1}-u_{1})$, $h_{2}=h\cdot (v_{2}-u_{2})$. But 
   then, Prop. \ref{firstorder2} implies 
   \begin{align*}
	h_{1}D_{1}(\tilde{f}_{t};\underline{u})(\underline{v})&+
	h_{2}D_{2}(\tilde{f}_{t};\underline{u})(\underline{v})\\
	&=h(D_{1}(\tilde{f}_{t};\underline{u})(\underline{v})(v_{1}-u_{1})+
	D_{2}(\tilde{f}_{t};\underline{u})(\underline{v})(v_{2}-u_{2}))
	=h\cdot f_{0}(u_{1},u_{2})
    \end{align*}
    proving the injectivity of the map.      
\end{proof}

\section{Lifting factorisations to $k\langle 
x,y\rangle$}\label{factorisation}

There is no general unique factorisation theorem in $k\langle 
x,y\rangle$, e.g. the two factorisations
\begin{displaymath}
    (xy+1)x=x(yx+1)
\end{displaymath}
are different. Moreover, a factorisation $f=gh$ in $k[x,y]$ of 
an element $f\in k\langle x,y\rangle$ cannot always be lifted back 
to $k\langle x,y\rangle$.
\par
Let $A=k\langle x,y\rangle/(f)$ and $A_{0}=k[x,y]/(f_{0})$ the natural 
commutativisation. Put $C_{f}:=Simp_{1}(A)$. 
\begin{Prop}
    Let $f\in k\langle x,y\rangle$ be reducible with factorisation 
    $f=gh$. Then for all $p_{1}\in C_{h}$ and $p_{2}\in C_{g}$ we have 
    \begin{displaymath}
	Ext_{A}^1(k(p_{1}),k(p_{2}))\neq 0
    \end{displaymath}
    where $A=k\langle x,y\rangle/(f)$.
\end{Prop}
\begin{proof}
    Let $p_{1}=(\alpha_1,\beta_1)$ and $p_{2}=(\alpha_2,\beta_2)$. Then 
    \begin{displaymath}
	D_{i}(f;p_{1})(p_{2})
	=g(p_{2})D_{i}(f;p_{1})(p_{2})+D_{i}(g;p_{1})(p_{2})h(p_{1})=0
    \end{displaymath}
    for $i=1,2$. As shown in section \ref{Jacobi} this implies that 
    $Ext_{A}^1(k(p_{1}),k(p_{2}))\neq 0$.
\end{proof}
We are interested in some sort of converse to this statement. Can we 
from geometrical data connected to the curve $C_{f}$ decide whether or 
not
there exists a factorisation $f=gh\in k\langle x,y\rangle$? We shall 
not try to give a complete solution to this problem, but in the following 
theorem we give a partial answer.   
\begin{Thm}
    Let $f\in k\langle x,y\rangle$ be such that $f_{0}$ is reduced and 
    has a proper factorisation $f_{0}=gh$ in $k[x,y]$. Assume
    \begin{displaymath}
	Ext_{A}^1(k(p_{1}),k(p_{2}))\neq 0
    \end{displaymath}
    for all $p_{1}\in C_{h}$, $p_{2}\in C_{g}$ and where $A=k\langle 
    x,y\rangle/(f)$. Then there exist non-commutative models
    $\hat{g}$ and $\hat{h}$ of $g$ resp. $h$, such that 
    \begin{displaymath}
	f-\hat{g}\hat{h}\in ([x,y])^2
    \end{displaymath}
\end{Thm}
\begin{proof}
    Initially we choose non-commutative models $\hat{g}$ and $\hat{h}$ 
    for $g$ and $h$ such that the degree of the models do not exceed 
    the degree of the original elements.
    Then we have 
    \begin{displaymath}
	f-\hat{g}\hat{h}=\sum_{i=1}^l g_{i}[x,y]h_{i}
    \end{displaymath}
    for some $g_{i},h_{i}\in k\langle x,y\rangle$ and such that
    $deg(g_{i})+deg(h_{i})<deg(f)-1$ for all $i=1,\ldots,l$.        
    For $p_{1}=(\alpha_1,\beta_1)\in C_{h}$ and $p_{2}
    =(\alpha_2,\beta_2)\in C_{g}$ we have 
    \begin{align*}
	D_{1}(f;p_{1})(p_{2})
	&=
	\sum g_{i}(\alpha_2,\beta_2)(\beta_{1}-\beta_{2})h_{i}(\alpha_1,\beta_1)\\	
	D_{2}(f;p_{1})(p_{2})
	&=
	\sum g_{i}(\alpha_2,\beta_2)(\alpha_{2}-\alpha_{1})h_{i}(\alpha_1,\beta_1)\\
    \end{align*}
    By assumption $D_{1}(f;p_{1})(p_{2})=D_{2}(f;p_{1})(p_{2})=0$. Since 
    $C_{g}$ and $C_{h}$ are non-empty curves with no common component 
    this implies that
    \begin{displaymath}
	(*)\qquad\sum_{i=1}^l g_{i}(\alpha_2,\beta_2)h_{i}(\alpha_1,\beta_1)=0, 
	\qquad (\alpha_2,\beta_2)\in C_{g}, 
	(\alpha_1,\beta_1)\in C_{h}
    \end{displaymath}
    Let us first treat the case $l=1$. Then we have 
    $g(\alpha_2,\beta_2)h(\alpha_1,\beta_1)=0$ for all 
    $(\alpha_2,\beta_2)\in C_{g}, (\alpha_1,\beta_1)\in C_{h}$ and  
    $deg(h_{1})<deg(h)$ or $deg(g_{1})<deg(g)$. 
    If $deg(g_{1})<deg(g)$ there exists 
    $(\alpha_2,\beta_2)\in C_{g}$ such that 
    $g_{1}(\alpha_2,\beta_2)\neq 0$ since $f$ is reduced. Hence $h_{1}(\alpha_1,\beta_1)=0$ 
    for all $(\alpha_1,\beta_1)\in C_{h}$ and $h_{1}-k_{1}\hat{h}\in 
    ([x,y])$ for some $k_{1}\in k\langle x,y\rangle$ and
    \begin{displaymath}
	f=\hat{g}\hat{h}+g_{1}[x,y]k_{1}\hat{h}+([x,y])^2
    \end{displaymath}
    hence the claim holds, since $\hat{g}+\hat{g}_{1}[x,y]\hat{k}_{1}$ 
    is a non-commutative model for $g$.
    \par
    Assume
    \begin{displaymath}
	\sum_{i=1}^l g_{i}(\alpha_2,\beta_2)h_{i}(\alpha_1,\beta_1)=0
    \end{displaymath}
    for all $(\alpha_1,\beta_1)\in C_{h}$ and for all $(\alpha_2,\beta_2)\in C_{g}$. 
    By symmetry we can assume $deg (h_{1})< deg (h)$ and 
    choose $(\alpha_{0},\beta_{0})$ such that $h_{1}(\alpha_{0},\beta_{0})\neq 0$. 
    Then for all $\alpha_{2}$ and $\beta_{2}$ 
    \begin{displaymath}
	g_{1}(\alpha_2,\beta_2)=\sum_{i=2}^l\lambda_{i}g_{i}(\alpha_2,\beta_2)
    \end{displaymath}
    for some $\lambda_{i}\in k$. Inserting this in (*) we get 
    \begin{displaymath}
	\sum_{i=2}^l(h_{i}(\alpha_1,\beta_1)+\lambda_{i}h_{1}(\alpha_1,\beta_1))g_{i}(\alpha_2,\beta_2)=0
    \end{displaymath}
    for all $\alpha_{2}$ and $\beta_{2}$. If we have started with ``$h_{1}$'' of minimal 
    degree we get
    \begin{displaymath}
	deg (h_{i}+\lambda_{i}h_{1})\le deg( h_{i})
    \end{displaymath}
    This procedure will of course lead us to the case $l=1$ which was 
    treated above.
\end{proof}

\section{The completion problem}

The main result of this section gives a partial
answer to the following problem, which indecomposable modules deform 
into simple modules? 

\subsection{Graphs with no complete cycles}\label{graphs without 
complete cycles}

\begin{Defn}
    Let $A$ be a $k$-algebra and $\V=\{V_{1},\ldots,V_{r}\}$ a finite set of simple 
    $A$-modules. The directed extension graph $Q\V$ of $\V$ has the modules 
    $V_{i}$ as its vertices and there is an arrow from $V_{i}$ to $V_{j}$ if and 
    only if $Ext_{A}^1(V_{i},V_{j})\neq 0$.
\end{Defn}
By a cycle in the extension graph $Q\V$ we shall mean a path, 
starting and 
ending at the same vertex. A cycle is said to be complete if it 
contains all the vertices of $Q\V$. The successor set
$E(\nu)$ of a vertex $\nu$ is the set of vertices $\omega$ of $Q\V$ such 
that there is a path from $\nu$ to $\omega$, and the precursor 
set $P(\omega)$ consists of all 
vertices $\nu\in\V$ such that there is a path from $\nu$ to $\omega$.
Notice that we always include $\nu$ both in the successor and the 
precursor set of $\nu$ and that $E(\nu)=\V$ for all 
vertices $\nu\in \V$ if and only if there 
exists a complete cycle in $\V$. 
\begin{Lem}\label{splitting}
    Let $Q\V$ be an extension graph with no complete cycles. Then 
    there exists a disjoint union $\V=\M\cup \N$ such that
    there are no arrows (in $Q\V$) from $Q\M$ to $Q\N$.
\end{Lem}
\begin{proof}
    Since there are no complete cycles in $Q\V$ there exists a vertex $\nu$ 
    such that $E(\nu)\neq \V$. Let $\M=E(\nu)$ and $\N=\V-\M$. Then 
    $\V=\M\cup\N$ is a disjoint union and since there are no arrows 
    out of a successor set, there are no arrows from $Q\M$ to $Q\N$.
\end{proof}

\subsection{Homomorphisms and extensions of indecomposable modules}
\label{homomorphisms}

Let $E$ be a non-simple, indecomposable $A$-module with a finite 
composition series. Let the cofiltration
\begin{displaymath}
    E=E_{r}\rightarrow E_{r-1}\rightarrow 
    \ldots\rightarrow E_{1}\rightarrow E_{0}=0
\end{displaymath}
be a cocomposition series of $E$ and let $V_{i}=ker\{E_{i}\rightarrow 
E_{i-1}\}$, $i=1,2,\ldots,r$ be the simple kernels. 
The set $\V=\{V_{1},\ldots,V_{r}\}=:Supp(E)$ is called the support of $E$. 
It only depends on $E$ and not on the 
cofiltration. Let $\Gamma$ be an ordering of the elements of $\V$, 
\begin{displaymath}
    \Gamma=(V_{i_{1}},V_{i_{2}},\ldots,V_{i_{r}})
\end{displaymath}
    An $A$-module $F$ with a finite cofiltration
\begin{displaymath}
    F=F_{r}\rightarrow F_{r-1}\rightarrow 
    \ldots\rightarrow F_{1}\rightarrow F_{0}=0
\end{displaymath}
such that $ker\{F_{j}\rightarrow F_{j-1}\}=V_{i_{j}}$, $j=1,2,\ldots,r$
is called an \textbf{iterated extension} 
of $Supp(E)$ defined by the ordering $\Gamma$. 
The set of all iterated 
$A$-module extensions defined by $\Gamma$ is denoted 
$Ind_{\Gamma}(A)$. It 
has an affine non-commutative scheme structure (see \cite{La4}).
\begin{Lem}\label{hom of disjoint}
    Let $M,N$ be finite dimensional $A$-modules with disjoint support.
    Then 
    \begin{displaymath}
	Hom_{A}(M,N)=0
    \end{displaymath}
\end{Lem}
\begin{proof}
    The proof is by induction on the cardinality of 
    $Supp(N)$, resp. $Supp(M)$.
    \par
    Assume first that $M$ and $N$ are simple modules.
    By Schurs lemma 
    \begin{displaymath}
	Hom_{A}(M,N)=0
    \end{displaymath}
    since $M$ and $N$ are non-isomorphic.
    Let $M$ be simple and $V\subseteq 
    N$. Then there is a short-exact sequence 
    \begin{displaymath}
	0\rightarrow V\longrightarrow N\longrightarrow N/V\rightarrow 0
    \end{displaymath}
    inducing a long-exact sequence 
    \begin{displaymath}
	0\rightarrow Hom_{A}(M,V)\rightarrow Hom_{A}(M,N)
	\rightarrow Hom_{A}(M,N/V)\rightarrow Ext^1_{A}(M,V) 
	\rightarrow \ldots
    \end{displaymath}
    By induction hypothesis $Hom_{A}(M,V)=Hom_{A}(M,N/V)=0$ and 
    therefore also $Hom_{A}(M,N)=0$. 
    \par
    A similar argument gives the induction step for 
    $Hom_{A}(\textrm{-},N)$.
\end{proof}
\begin{Lem}\label{endomorphisms}
    Let $E$ be a non-simple, indecomposable $A$-module given by a 
    short-exact sequence
    \begin{displaymath}
	0\rightarrow M\longrightarrow E\longrightarrow N\rightarrow 0
    \end{displaymath}
    with $Supp(M)\cap Supp(N)=\emptyset$. Suppose $End_{A}(M)\simeq 
    End_{A}(N)\simeq k$. Then $End_{A}(E)\simeq k$.
\end{Lem}
\begin{proof}
    We have the following diagram of $k$-vector spaces
    \begin{displaymath}
	\diagram
	     &0\dto               &0\dto               &0\dto\\
	0\rto&Hom_{A}(N,M)\dto\rto&Hom_{A}(N,E)\dto\rto&Hom_{A}(N,N)\dto\\
	0\rto&Hom_{A}(E,M)\dto\rto&Hom_{A}(E,E)\dto\rto&Hom_{A}(E,N)\dto\\
	0\rto&Hom_{A}(M,M)    \rto&Hom_{A}(M,E)    \rto&Hom_{A}(M,N)\\
	\enddiagram
    \end{displaymath}
    By Lemma \ref{hom of disjoint} 
    \begin{displaymath}
	Hom_{A}(N,M)=Hom_{A}(M,N)=0
    \end{displaymath} 
    and there are inclusions
    \begin{align*}
	&0\rightarrow Hom_{A}(N,E)\longrightarrow Hom_{A}(N,N)\simeq k\\
	&0\rightarrow Hom_{A}(E,M)\longrightarrow Hom_{A}(M,M)\simeq k\\
    \end{align*}
    Since the given short-exact sequence is not split, this implies 
    that
    \begin{displaymath}
	Hom_{A}(N,E)=Hom_{A}(E,M)=0
    \end{displaymath}
    Thus $End_{A}(E)\simeq k$.
\end{proof}
We have the following characterisation of 
indecomposable modules.
\begin{Prop}
    Let the $A$-module $E$ be an iterated extension such that $Supp(E)$ 
    contains no multiple 
    simple modules. Then $E$ is indecomposable if and only if 
    $End_{A}(E)\simeq k$.
\end{Prop}
\begin{proof}
    The only if part follows from Lemma \ref{endomorphisms} and a 
    straightforward induction argument. If $E$ is decomposable we can 
    write $E\simeq E^{\prime}\oplus E^{\prime\prime}$ as a proper 
    sum, proving that $dim_{k}End_{A}(E)\ge 2$.
\end{proof}
In addition to these general facts about homomorphisms we also need 
some basic results about infinitesimal homomorphisms, i.e. about the various 
ext-groups.
\begin{Lem}\label{vanishing ext}
    Let $\V=\{V_{1},\ldots,V_{r}\}$ be a family of simple $A$-modules and 
    let $Q\V$ be the corresponding extension graph. 
    Suppose there are no complete cycles in $Q\V$, and choose $\nu\in Q\V$ 
    such that $E(\nu)\neq \V$. If $M$  
    an $A$-module with $Supp(M)=E(\nu)$ and $N$ 
    an $A$-module such that $Supp(N)\cap E(\nu)=\emptyset$, then 
    \begin{displaymath}
	Ext^1_{A}(M,N)=0
    \end{displaymath}
\end{Lem}
\begin{proof}
    Let 
    \begin{displaymath}
	0=M_{0}\subset M_{1}\subset M_{2}\subset \ldots\subset M_{m}=M
    \end{displaymath}
    be a filtration of $M$ such that
    \begin{displaymath}
	M_{i}/M_{i-1}\simeq V_{j_{i}}\,, \quad i=1,2,\ldots,m
    \end{displaymath}
    Let 
    \begin{displaymath}
	0=N_{0}\subset N_{1}\subset N_{2}\subset \ldots\subset N_{n}=N
    \end{displaymath}
    be a similar filtration of $N$. We claim that 
    \begin{displaymath}
	Ext^1_{A}(M_{i},N_{k})=0 \quad i=1,2,\ldots,m\quad k=1,2,\ldots,n
    \end{displaymath}
    In fact there are long-exact sequences
    \begin{align*}
	0\rightarrow &Hom_{A}(V_{j_{i}},N_{k})\longrightarrow 
	Hom_{A}(M_{i},N_{k})\longrightarrow Hom_{A}(M_{i-1},N_{k})\rightarrow\\
	& Ext^1_{A}(V_{j_{i}},N_{k})\longrightarrow 
	Ext^1_{A}(M_{i},N_{k})\longrightarrow Ext^1_{A}(M_{i-1},N_{k})\\
    \end{align*}
    Now by Lemma \ref{hom of disjoint} $Hom_{A}(M_{i-1},N_{k})=0$ since 
    $Supp(M_{i-1})\cap 
    Supp(N_{k})=\emptyset$. For any $V_{s}\in Supp(M)$ and $V_{t}\in 
    Supp (N)$ we have by assumption $Ext^1_{A}(V_{s},V_{t})=0$ and an 
    induction argument using the long-exact sequence above and
    similar sequences for the extensions in $N$ gives the result.
\end{proof}

\subsection{A completion theorem}\label{completion}

Let $\V=\{V_{1},\ldots,V_{r}\}$ be a finite family of simple $A$-modules 
such that the extension graph $Q\V$ has no complete cycles. By Lemma 
\ref{splitting} there is a disjoint union 
$\V=\M\cup \N$ with no arrows in $Q\V$ from $Q\M$ to $Q\N$. 
\begin{Defn}
    Let $\V$ be a finite family of simple $A$-modules 
    and let $\V=\M\cup \N$ be a disjoint 
    union with no arrows in $Q\V$ from $Q\M$ to $Q\N$.
    An ordering $\Gamma=(V_{j_{1}},\ldots,V_{j_{r}})$ of $\V$ is 
    said to respect the union $\V=\M\cup \N$ if for some integer $l$ 
    $V_{j_{1}},\ldots,V_{j_{l}}\in \N$ and 
    $V_{j_{l+1}},\ldots,V_{j_{r}}\in \M$.
\end{Defn}
The set of indecomposable $A$-modules defined as iterated extensions 
of the simple modules $V_{1},\ldots,V_{r}$, in accordance with the 
ordering $\Gamma$, is denoted $Ind_{\Gamma}(A)$ (see \cite{La4} for 
details).
\par
The next proposition is crucial for the Completion Theorem \ref{main 
theorem}. It will be proved in section \ref{moving simple modules}.
\begin{Prop}\label{exists extension}
    Let $E$ be an indecomposable $A$-module with $Supp(E)=\V$ and 
    $End_{A}(E)\simeq k$ and suppose there are 
    no complete cycles in $Q\V$. Then there exists a short-exact 
    sequence
    \begin{displaymath}
	0\rightarrow M\longrightarrow E\longrightarrow N\rightarrow 0
    \end{displaymath}
    such that there are no arrows in $Q\V$ from $Q\M$ to $Q\N$, where
    $\M=Supp(M)$ and $\N=Supp(N)$.
\end{Prop}
By choosing any cofiltration of $M$ and $N$ we see that
there is an equivalent conclusion for this proposition, saying that
with the same assumptions, there exists a $\Gamma$ which 
respects the union $\V=\M\cup \N$ and such that $E\in Ind_{\Gamma}(A)$. 
\par
Let $\Gamma$ be some ordering of $\V$. Denote by 
\begin{displaymath}
    Simp_{\Gamma}(A)=Ind_{\Gamma}(A)\cup Simp_{n}(A)
\end{displaymath}
where $n$ is the $k$-dimension of the $A$-modules in $Ind_{\Gamma}(A)$. 
Laudal has shown (see \cite{La4}) that $Simp_{\Gamma}(A)$ has a 
non-commutative scheme structure, called a completion of $Simp_{n}(A)$, 
represented by a sheaf of $k$-algebras. There are of course several 
completions of $Simp_{n}(A)$, but we do not know whether or not there exists a 
non-commutative scheme which at the same time parametrises all possible completions 
of $Simp_{n}(A)$. Nevertheless, for any indecomposable $A$-module 
with $End_{A}(E)\simeq k$ there exists a formal moduli $H(E)$ for the 
deformation functor $Def_{E}$.
\par
For a short-exact sequence
\begin{displaymath}
    \xi:\qquad 0\rightarrow M\longrightarrow E\longrightarrow N\rightarrow 0
\end{displaymath}
we consider the functor $Def_{\xi}$, containg all 
liftings of $E$ which preserves the extension structure, i.e. liftings 
of $E$ which are extensions of liftings of $N$ by liftings of $M$. A 
more precise description can be found in section \ref{simple 
deformations}. The 
formal moduli for $Def_{\xi}$ is denoted $H(\xi)$ and there is a 
map of formal $k$-algebras
\begin{displaymath}
    \tilde{\xi}:H(E)\longrightarrow H(\xi)
\end{displaymath}
The kernel of this map is our candidate for deformations pointing into 
the simple locus. But even if $\tilde{\xi}$ induces an injection on
tangent level, $E$ may deform into simple modules. To avoid such 
``singularities'' in the completion of $Simp_{n}(A)$ at $E$ we make 
the following definition.
\begin{Defn}
    An indecomposable $A$-module $E$ with $End_{A}(E)\simeq k$ is 
    called a smooth completion point of $Simp_{n}(A)$ if for 
    allnon-split exact sequences
    \begin{displaymath}
	\xi:\qquad 0\rightarrow M\longrightarrow E\longrightarrow N\rightarrow 0
    \end{displaymath}
    the map $\tilde{\xi}$ is not a monomorphism on the tangent 
    level.
\end{Defn}
The last ingredient of the completion theorem is the bridge between 
deformation theory and the homological properties of the extension 
graph $Q\V$. It is builted in the next proposition. The proof is 
postponed to section \ref{simple deformations}.
\begin{Prop}\label{non-vanishing ext}
    Suppose the indecomposable $A$-module $E$ is a smooth completion 
    point of $Simp_{n}(A)$, and let 
    \begin{displaymath}
	\xi:\qquad 0\rightarrow M\longrightarrow E\longrightarrow N\rightarrow 0
    \end{displaymath}
    be some exact sequence presenting $E$. Then
    \begin{displaymath}
	Ext_{A}^1(M,N)\neq 0
    \end{displaymath}
\end{Prop}
The main result of this sections is the completion theorem, giving a 
criterion for when the indecomposable $A$-module $E$ deforms into 
simple modules.
\begin{Thm}\label{main theorem}
    Let $A$ be an associative $k$-algebra and $E$ an indecomposable $A$-module, 
    such that $End_{A}(E)\simeq k$. Let 
    $\V=Supp(E)=\{V_{1},\ldots,V_{r}\}$ be the associated simple 
    modules and let $Q\V$ be the extension graph. Suppose $E$ is a 
    smooth completion point of $Simp_{n}(A)$. Then there exist a complete 
    cycle in $Q\V$.
\end{Thm}
\begin{proof}
    Assume for a contradiction that there are no complete cycles in 
    $Q\V$. Then by Lemma \ref{splitting} there is a disjoint 
    union $\V=\M\cup\N$ with no arrows in $Q\V$ from $Q\M$ to $Q\N$.  
    Then there exist a short-exact sequence  
    \begin{displaymath}
	\xi:\qquad 0\rightarrow M\longrightarrow E\longrightarrow N\rightarrow 0
    \end{displaymath}
    There are no arrows in $Q\V$ from $Q\M$ to $Q\N$ and by 
    Proposition \ref{vanishing ext} this implies that
    \begin{displaymath}
	Ext_{A}^1(M,N)=0
    \end{displaymath}
    contradicting the assertion of Proposition \ref{non-vanishing ext}. 
\end{proof}
Notice that the existence of a complete cycle is not a sufficient 
condition for existence of simple modules of the given dimension. 
Consider e.g. the plane curve given by the equation 
$f=x^2-1+[x,y]^2y$. It has complete 2-cycles, but there are no 2-dimensional 
simple modules. In fact it is easy to construct examples of two 
$k$-algebras $A$ and $A^{\prime}$ such that $Simp_{1}(A)\simeq 
Simp_{1}(A^{\prime})$ and with isomorphic extension graphs, but such 
that the higher dimensional simple structures are completely different. 

\subsection{Proof of Proposition \ref{non-vanishing ext}}
\label{simple deformations}

Let $N$ and $M$ be left $A$-modules and consider an extension $E$ of 
$N$ 
by $M$, given by the short-exact sequence
\begin{displaymath}
    \xi:\qquad 0\rightarrow M\buildrel\iota\over\longrightarrow 
    E\buildrel\pi\over\longrightarrow 
    N\rightarrow 0
\end{displaymath}
Let $\underline{l}$ denote the category of commutative local artinian 
$k$-algebras with residue field $k$. We associate to the exact 
sequence $\xi$ a 
deformation functor 
\begin{displaymath}
    Def_{\xi}:\underline{l}\rightarrow \underline{Sets}
\end{displaymath}
defined as follows. For any pointed local artinian $k$-algebra $\rho: 
R\rightarrow k$ the set $Def_{\xi}(R)$ consists of isomorphism 
classes of short-exact 
sequences of left $A\otimes R^{op}$-modules
\begin{displaymath}
    \xi_{R}:\qquad 0\rightarrow M\otimes R\buildrel\iota_{R}\over\longrightarrow 
    E\otimes R\buildrel\pi_{R}\over\longrightarrow 
    N\otimes R\rightarrow 0
\end{displaymath}
such that the $A$-module structure of $E\otimes R$, given by a 
$k$-algebra homomorphism 
\begin{displaymath}
    \eta_{E}^{\prime}:A\longrightarrow End_{R}(E\otimes R)
\end{displaymath}
is compatible with some $A$-module structure of $N\otimes R$ and 
$M\otimes R$, and such that $\xi_{R}\otimes_{R}k\simeq \xi$. 
\par
There is a forgetful morphism of deformation functors
\begin{displaymath}
    Def_{\xi}\longrightarrow Def_{E}
\end{displaymath}
inducing a $k$-linear map of tangent spaces
\begin{displaymath}
    Def_{\xi}(k[\epsilon])\longrightarrow Def_{E}(k[\epsilon])\simeq 
    Ext_{A}^1(E,E)
\end{displaymath}
The image of the map is denoted 
$Ext_{A}^1(E,E)_{\xi}\subset Ext_{A}^1(E,E)$. Our concern is to 
determine when the map is surjective.
\par
An element $\alpha\in Def_{\xi}(k[\epsilon])$ is
given by a triple $\alpha=(\beta,\eta,\mu)$, 
where 
\begin{displaymath}
    \beta\in Ext_{A}^1(E,E)\quad\eta\in Ext_{A}^1(M,M)\quad\mu\in Ext_{A}^1(N,N)
\end{displaymath}
and such that $\iota_{\ast}(\eta)=\iota^{\ast}(\beta)$ and
$\pi^{\ast}(\mu)=\pi_{\ast}(\beta)$.
\par
\begin{displaymath}
    \diagram
    Ext^1_{A}(M,M)\dto^{\iota_{\ast}}&\dlto^{\iota^{\ast}} 
    Ext^1_{A}(E,E)\ddto^{\chi}\drto_{\pi_{\ast}}&\dto^{\pi^{\ast}} 
    Ext^1_{A}(N,N)\\
    Ext^1_{A}(M,E)\drto_{\pi_{\ast}}&
    &\dlto^{\iota^{\ast}} Ext^1_{A}(E,N)\\
    &Ext^1_{A}(M,N)&\\
    \enddiagram
\end{displaymath}
Refering to the diagram we put 
$\chi:=\iota^{\ast}\circ\pi_{\ast}=\pi_{\ast}\circ\iota^{\ast}$.
\begin{Lem}
    If $Ext_{A}^1(E,E)_{\xi}\neq Ext_{A}^1(E,E)$, then 
    $Ext^1_{A}(M,N)\neq 0$.
\end{Lem}
\begin{proof}
    Assume $Ext^1_{A}(M,N)= 0$. Then the above diagram 
    reduces to
    \begin{displaymath}
	\diagram
	Ext^1_{A}(M,M)\dto^{\iota_{\ast}}&\dlto^{\iota^{\ast}} 
	Ext^1_{A}(E,E)\drto_{\pi_{\ast}}&\dto^{\pi^{\ast}} 
	Ext^1_{A}(N,N)\\
	Ext^1_{A}(M,E)& & Ext^1_{A}(E,N)\\
	\enddiagram
    \end{displaymath}
    where $\iota_{\ast}$ and $\pi^{\ast}$ are surjections. Thus for 
    any $\beta\in Ext_{A}^1(E,E)$ there 
    exist $\eta\in Ext^1_{A}(M,M)$ and $\mu\in Ext^1_{A}(N,N)$ 
    such that $\iota_{\ast}(\eta)=\iota^{\ast}(\beta)$ and
    $\pi^{\ast}(\mu)=\pi_{\ast}(\beta)$. Consequently 
    $Ext_{A}^1(E,E)_{\xi}= Ext_{A}^1(E,E)$.
\end{proof}
\noindent\textit{Proof of Proposition \ref{non-vanishing ext}.}
    The condition of the lemma says that the map 
    \begin{displaymath}
	\tilde{\xi}:H(E)\longrightarrow H(\xi)
    \end{displaymath}
    is not a momomorphism on tangent level, i.e. $E$ is a smooth 
    completion point of $Simp_{n}(A)$. The conclusion is the 
    non-vanishing of $Ext_{A}^1(M,N)$.
\begin{flushright}
    \textsquare
\end{flushright}

\subsection{Proof of Proposition \ref{exists extension}}
\label{moving simple modules}

In this section we shall give a proof of Proposition \ref{exists 
extension}. We start with the so-called Moving Lemma. It says that when 
certain ext-group vanishes, we can ``move'' components of the 
cofiltration of a given indecomposable module.
\begin{Lem}\label{moving lemma}
    Consider the diagram of short-exact sequences of $A$-modules
    \begin{displaymath}
	\diagram
	      &                &       & 0 \dto                   & & \\
	      &                &       & N \dto                   & & \\
	0\rto & M^{\prime}\rto & M\rto & M^{\prime\prime}\rto\dto & 0 
	& \\
	      & 0\rto          & W\rto & N^{\prime}\rto\dto       & 
	W^{\prime}\rto & 0 \\
	      &                &       & 0                        & &  \\
	\enddiagram
    \end{displaymath}
    Suppose $Ext_{A}^1(W,N)=0$. Then there exists a similar diagram 
    with $W$ and $N$ interchanged:
    \begin{displaymath}
	\diagram
	      &                &       & 0 \dto                   & & \\
	      &                &       & W \dto                   & & \\
	0\rto & M^{\prime}\rto & M\rto & M^{\prime\prime}\rto\dto & & \\
	      & 0\rto          & N\rto & N^{\prime\prime}\rto\dto & 
	W^{\prime}\rto & 0 \\
	      &                &       & 0                        & &  \\
	\enddiagram
    \end{displaymath}
\end{Lem}
\begin{proof}
    Let $L_{1}$ be the kernel of the composed map 
    \begin{displaymath}
	M^{\prime\prime}\longrightarrow N^{\prime}\longrightarrow W^{\prime}
    \end{displaymath}
    A simple diagram chasing argument shows that $L_{1}$ 
    is an extension of $W$ by 
    $N$, i.e. there is a short-exact sequence
    \begin{displaymath}
	0\rightarrow N\longrightarrow L_{1}\longrightarrow 
	W\rightarrow 0
    \end{displaymath}
    By assumption this sequence splits and there is a section $s:W\rightarrow L_{1}$ which composed 
    with the surjection $L_{1}\rightarrow W$ gives the identity on $W$.
    The composed map
    \begin{displaymath}
	W\buildrel s\over\longrightarrow L_{1}\longrightarrow M^{\prime\prime}
    \end{displaymath}
    is clearly injective. Let $N^{\prime\prime}$ be the cokernel.
    The composed map
    \begin{displaymath}
	N\longrightarrow M^{\prime\prime}\longrightarrow N^{\prime\prime}
    \end{displaymath}
    is again injective, with cokernel 
    $W^{\prime}$ and we are done.
\end{proof}
Let $\V=\{V_{1},\ldots,V_{r}\}$ be any set of simple $A$-modules with 
extension graph $Q\V$. Let $E$ be any $A$-module with $Supp(E)=\V$. 
Let $\V_{0}\subset \V$ be a subfamily such that there are no arrows 
in $\V$ from any 
element of $\V_{0}$ to any element outside $\V_{0}$. i.e. if 
$V_{s}\in\V_{0}$ and $V_{t}\notin\V_{0}$, then 
$Ext_{A}^1(V_{s},V_{t})=0$. 
\begin{Lem}\label{submodule criterion}
    Let $E$ be an $A$-module with $Supp(E)=\V$. 
    Let $\V_{0}\subset \V$ be a subfamily as described above, with no 
    arrows out of $\V_{0}$. Then there exists a submodule 
    $E^{\prime}\subset E$ such that $Supp(E^{\prime})=\V_{0}$.
\end{Lem}
\begin{proof}
    Let
    \begin{displaymath}
	E=E_{r}\buildrel {f_{r}}\over\longrightarrow E_{r-1}
	\longrightarrow \ldots 
	\buildrel {f_{1}}\over\longrightarrow E_{1}
    \end{displaymath}
    be a cofiltration of $E$, where $f_{i}$ is surjective with kernel
    \begin{displaymath}
	ker(f_{i})=V_{i} \quad i=1,2,\ldots,r
    \end{displaymath}
    and where we put $K_{i}=ker\{E\rightarrow E_{i}\}$, $i=1,\ldots,r$.
    \par
    Suppose $V_{i-1}\in\V_{0}$ and $V_{i}\not\in\V_{0}$. Applying Lemma 
    \ref{moving lemma} to the diagram
    \begin{displaymath}
	\diagram
	      &           &             & 0 \dto          &             &   \\
	      &           &             & V_{i} \dto      &             &   \\
	0\rto & K_{i}\rto & E\rto       & E_{i}\rto\dto   &             &   \\
	      & 0\rto     & V_{i-1}\rto & E_{i-1}\rto\dto & E_{i-2}\rto & 0 \\
	      &           &             & 0               &             &   \\
	\enddiagram
    \end{displaymath}
    we obtain a new cofiltration corresponding to $\Gamma^{\prime}$ 
    derived from $\Gamma$ by interchanging $V_{i-1}$ and $V_{i}$. 
    Repeat this procedure until there are no members of $\V_{0}$ 
    proceedeing any non-members of $\V_{0}$. Let $\Gamma^{\prime\prime}$ 
    be the corresponding ordering. Then $E^{\prime}=ker\{E\rightarrow 
    E_{j}^{\prime\prime}\}$ for some specific $j$ will do the job.
\end{proof}
Notice that the same argument proves the ``dual'' lemma:
\begin{Lem}
    Let $E$ be a $A$-module with $Supp(E)=\V=\{V_{1},\ldots,V_{r}\}$. 
    Let $\V_{1}\subset \V$ be a subfamily, with no 
    arrows into $\V_{1}$, i.e. for any $V_{s}\notin\V_{1}$ and 
    $V_{t}\in\V_{1}$, we have $Ext_{A}^1(V_{s},V_{t})=0$. Then there 
    exists a quotient module 
    $E\rightarrow E^{\prime\prime}$ such that $Supp(E^{\prime\prime})=\V_{1}$.
\end{Lem}    
\begin{proof}
    A dual, but quite similar argument as in the proof of lemma 
    \ref{submodule criterion}.
\end{proof}
\noindent\textit{Proof of Proposition \ref{exists extension}.}
    Let 
    \begin{displaymath}
	0\rightarrow M\longrightarrow E\longrightarrow N\rightarrow 0
    \end{displaymath}
    be some extension defining $E$. If there are no arrows in $Q\V$ 
    from $Q\M$ to $Q\N$ we are done.
    \par
    Assume in contrary that there exist an arrow in $Q\V$ from $Q\M$ to $Q\N$. 
    Let $\nu\in Q\N$ be the target of this arrow and consider the successor 
    set $E(\nu)$. By Lemma \ref{submodule criterion} there is a 
    submodule $N^{\prime}\hookrightarrow N$ with 
    $Supp(N^{\prime})=E(\nu)$ and a diagram
    \begin{displaymath}
	\diagram
	      &       &       & 0 \dto           &   \\
	      &       &       & N^{\prime} \dto  &   \\
	0\rto & M\rto & E\rto & N\rto\dto        & 0 \\
	      &       &       & N/N^{\prime}\dto &   \\
	      &       &       & 0                &   \\
	\enddiagram
    \end{displaymath}
    The kernel of the composed map
    \begin{displaymath}
	K=ker \{E\longrightarrow N\longrightarrow N/N^{\prime}\}
    \end{displaymath}
    is an extension of $N^{\prime}$ by $M$
    \begin{displaymath}
	0\rightarrow M\longrightarrow K\longrightarrow N^{\prime}\rightarrow 0
    \end{displaymath}
    and there is of course an exact sequence 
    \begin{displaymath}
	0\rightarrow K\longrightarrow E\longrightarrow Q\rightarrow 0
    \end{displaymath}
    where $Q=N/N^{\prime}$ is the quotient module.
    By construction 
    \begin{displaymath}
	Supp(K)=Supp(M)\cap E(\nu)
    \end{displaymath}
    \par
    We have ``moved'' the successor graph $E(\nu)$ from the epi side to the 
    mono side of certain extension sequences. Repeat this procedure, i.e. 
    construct new extension sequences
    \begin{displaymath}
	0\rightarrow K_{j}\longrightarrow E\longrightarrow Q_{j}\rightarrow 0
	\quad j=1,2,\ldots,l
    \end{displaymath}
    until there are no more arrows in the extension graph from the mono side 
    to the epi side. The precise meaning of the last statement is that 
    for any $V_{s}\in Supp(K_{l})$ and $V_{t}\in Supp(Q_{l})$, we have 
    $Ext_{A}^1(V_{s},V_{t})=0$.
    \par 
    There is a dual procedure to the one described, moving a precursor 
    subgraph $P(\omega)$ from 
    the mono side to the epi side.
    \par
    The procedures work as long as $E(\nu)\neq \N$ or $P(\omega)\neq 
    \M$. If $E(\nu)=\N$ and $P(\omega)=\M$ the procedures would in fact 
    collapse the exact sequence. So suppose $E(\nu)=\N$ and $P(\omega)= 
    \M$ for all $\nu\in Q\N$ and $\omega\in Q\M$. 
    There are two possibilities: 
    \par
    1) There are only 
    arrows from mono to epi side, i.e. for any $V_{s}\in Supp(N)$ and 
    $V_{t}\in Supp(M)$, we have 
    $Ext_{A}^1(V_{s},V_{t})=0$. By lemma \ref{vanishing ext} this implies 
    that $Ext_{A}^1(N,M)=0$ and the sequence splits. Thus the 
    short-exact sequence
    \begin{displaymath}
	0\rightarrow N\longrightarrow E\longrightarrow M\rightarrow 0
    \end{displaymath}
    will do the job.
    \par
    2) There are arrows in both directions. But then it is easily seen that 
    there exists a complete cycle, in contradiction to our assumptions. 
\begin{flushright}
    \textsquare
\end{flushright}


\section{The geometry of plane noncommutative curves}\label{conic}

To illustrate the theory of the previous sections we consider 
non-commutative models of plane 
quadrics and cubics, their finite-dimensional simple modules and the 
relations defined in the previous sections. 

\subsection{Classification of quadrics}

Up to a linear shift of basis a general plane noncommutative quadric
can be written
\begin{displaymath}
    A=k\langle x,y\rangle/(\lambda_{1} x^2+\lambda_{2} y^2+
    \delta [x,y]+ex+fy+g)
\end{displaymath}
with coefficients in $k$.
Going through all possibilities for the coefficients we achieve the 
following classification of plane noncommutative quadrics, up to 
isomorphism. 
\begin{itemize}
    \item[1)] The smooth case, $x^2+y^2-1+\delta[x,y]$. 
    This curve is isomorphic to the quantised Weyl algebra $xy-qyx-1$,
    whith $q=\frac{\delta +i}{\delta -i}$ for $\delta\neq i$.
    \item[2)] The singular case $x^2+y^2+\delta[x,y]$ which is better 
    known as the quantum plane $xy-qyx$ with $q\not=1$ and the same 
    relation between $q$ and $\delta$ as above.
    \item[3)] The degenerate cases, including two parallell lines 
    $x^2+ex+\delta[x,y]$, a double line ($e=0$) and a simple line 
    $x+\delta[x,y]$.
    \item[4)] The first Weyl algebra $1+\delta[x,y]$, $\delta\neq 0$.
    \item[5)] The affine plane $[x,y]$
\end{itemize}
The Weyl algebra has no finite-dimensional representations at all and 
the geometry of the commutative polynomial ring in two variables is 
well-known. The remaining cases are 1)-3) and we shall treat each case
separately.

\subsection{The quantum plane}\label{node}

For the quantum plane $xy=qyx$, $q\not= 1$, the following result holds. 
\begin{Prop}
    Let $A=k\langle x,y\rangle/(xy-qyx)$, where we assume $q\not=1$. 
    If $q$ is a primitive 
    $m$'th root of unity, the only simple modules of $A$ of dimension 
    strictly greater than one, are of dimension 
    $m$ and they correspond to 
    $(\lambda,\gamma)\in k^{\ast}\times k^{\ast}$, given by 
    \begin{displaymath}
	x\mapsto
	\bmatrix 
	\lambda\\ 
	& q\lambda\\ 
	&&\ddots\\ 
	&&& q^{m-1}\lambda
	\endbmatrix
	\quad
	y\mapsto
	\bmatrix 
	0&&&\gamma\\ 
	1&&&\\ 
	& \ddots&&\\ 
	& & 1&0
	\endbmatrix
    \end{displaymath}
    For other $q$ all finite dimensional simple representations are 
    1-dimensional.
\end{Prop}
\begin{proof}
    If $q\not=0$ the simple modules are well known (see e.g. \cite{DJ} 
    or \cite{J3}). For $q=0$ the defining relation is
    $$xy=0$$
    Now for a simple left $k\langle x,y\rangle/(xy)$-module $V$ the 
    submodule $yV$ is either 0 or $V$. Therefore
    $$
    Simp_1(A)=k\cup k\quad\text{and}\quad
    Simp_n(A)=\varnothing\;\;n>1\,
    $$
\end{proof}
The relation for this curve is the variety 
$\R\subset C\times C\subset \mathbf{A}^4$ given by the $k$-algebra 
$(k[x,y]/(f)\otimes k[u,v]/(f))/(qy-v,x-qu)$. The relation is 
of degree 1 and given by the map
\begin{displaymath}
    \Phi(u,v)=
    \left\{
    \begin{array}{ll}
	(0,q^{-1}v)& \textrm{if}\quad u=0\\
	(qu,0)& \textrm{if}\quad v=0\\
    \end{array}
    \right.
\end{displaymath} 
It is of finite order $m$ if and only if $q$ is a 
primitive $m$'th root of unity. In that case the completion of the simple locus
$Simp_{n}(A)$ is precisely the multi-extensions given by the 
finite cycles in the extension graphs. Thus 
$Simp_{1}(A)$ is singular and $Simp_{m}(A)$ is smooth (if it is 
non-empty).

\subsection{The smooth quadric}\label{smooth quadric}

For the smooth case we have a similar situation. The 
noncommutative \lq\lq conic section \rq\rq
\begin{displaymath}
    x^2+y^2-1+\delta[x,y]=0
\end{displaymath}
is isomorphic to the quantized Weyl algebra
\begin{displaymath}
    xy-qyx=1
\end{displaymath}
with $q=\frac{\delta-i}{\delta+i}$. Substituting 
$\tilde{y}=y-\frac{1}{1-q}x^{-1}$ we can transform the quantum plane 
into the quantised Weyl algebra (the cases where $x$ does not act as 
an automorphism are treated in \cite{J3}). Thus for $q$ an $m$'th root of unity we 
have simple modules of dimension $m$. In this case both $Simp_{1}(A)$ 
and $Simp_{m}(A)$ are smooth.
\par
For $\delta^2\neq -1$ the relation $\R$ is well-defined and of degree 1 
on $C$. It defines a linear map
\begin{displaymath}
    \Phi\matrV{u_{1}}{u_{2}}
    =\2matr{\frac{\delta^2-1}{\delta^2+1}}{\frac{-2\delta}{\delta^2+1}}
    {\frac{2\delta}{\delta^2+1}}{\frac{\delta^2-1}{\delta^2+1}}
    \matrV{u_{1}}{u_{2}}
\end{displaymath}
Suppose the ground field $k=\mathbf{C}$ is the complex numbers. Then 
for $\delta$ real the relation is precisely rotation by a fixed 
angle on a circle. 
For $\delta$ imaginary with $\delta\neq\pm i$ the 
relation is just hyperbolic rotation by a certain hyperbolic 
angle.
\par
Notice also that if $q$ is a primitive root of 1, 
$Ext_{A}^1(V,W)=0$ if $V$ or $W$ are non-isomorphic $m$-dimensional simple modules, 
by Proposition \ref{PI}. 
\par
As a brief remark we like to pay attention to the algebra given by the 
relation $xy+q[x,y]^2yx=1$ with $q\neq 0$. This is an other model for the smooth 
affine curve $xy-1$. In this case there are no non-trivial extensions 
between 1-dimensional representations.
\par
Let 
\begin{displaymath}
    \phi:k\langle x,y\rangle/(xy+q[x,y]^2yx-1)\rightarrow 
    M_{2}(\Gamma)
\end{displaymath}
be a 2-dimensional representation. By the Cayley-Hamilton 
theorem we can write 
\begin{displaymath}
    \phi(xy+q[x,y]^2yx-1)=XY-qd_{[X,Y]}YX-I
\end{displaymath}
where $\phi(x)=X$ and $\phi(y)=Y$. Now $X$ and $Y$ with the 
relation $XY-qd_{[X,Y]}YX-I=0$ generates $M_{2}(k(p))$ for $p\in 
Simp_{1}(\Gamma)$ if and only if 
$-qd_{[X,Y]}=1$. Using this fact and lemma \ref{traceformulas} we get
\begin{align*}
    XY-qd_{[X,Y]}YX-1
    &=XY-qd_{[X,Y]}(-XY+t_{X}Y+t_{Y}X+t_{XY}-t_{X}t_{Y})-1\\
    &=(1+qd_{[X,Y]})XY-qd_{[X,Y]}(t_{X}Y+t_{Y}X+t_{XY}-t_{X}t_{Y})-1
\end{align*}
and it follows that $qd_{[X,Y]}=-1$, $t_{X}=t_{Y}=0$ and $t_{XY}=1$. Substituting these 
values into the Formanek element we obtain the equality 
$4d_{X}d_{Y}=1-\frac{1}{q}$. Thus for $q\neq 1$ $Simp_{2}(R)$ is a smooth 
curve. In the singular case $q=1$ $Simp_{2}(R)$ consists of two lines 
with a normal crossing.

\subsection{The degenerate case}

For the simple line $x+\delta[x,y]$ we have \cite{DJ}
\begin{displaymath}
    Simp_n(A)=
    \begin{cases}
	\mathbf{A}^1& n=1\\
	\emptyset & n>1\\
    \end{cases}
\end{displaymath}
In the defining relation 
\begin{displaymath}
    x^2+\delta xy-\delta yx+ex=0, \qquad \delta\neq 0
\end{displaymath}
for the double line, $x$ is a normal element $(xA=Ax)$, and acts as $x=0$ or as
an automorphism on a simple module. The 1-dimensional simple modules correspond 
to $x=0$ or $x=-e$, i.e. 2 lines. If $x$ acts as an automorphism we write $z$ 
for $x^{-1}$ and the defining relation becomes
\begin{displaymath}
    1+\delta yz-\delta zy+ez=0
\end{displaymath}
This is either isomorphic to the Weyl algebra ($e=0$) or a simple 
line as described above ($e\not=0$). In both cases we have 
$Simp_n(A)=\emptyset$ for $n>1$.
\par
Notice that for the simple as well as for the double line, there exist non-trivial 
extensions of simple 1-dimensional representations:
If the defining relation is a simple line $x+\delta[x,y]$ the 
non-trivial extensions are given by ordered pairs
\begin{displaymath}
    [(0,\beta),(0,\beta + \frac{1}{\delta})]
\end{displaymath}
In the two-line case $x^2+ex+\delta[x,y]$ the non-trivial extensions are 
located in the ordered pairs
\begin{displaymath}
    [(0,\beta),(0,\beta + \frac{e}{\delta})]
\end{displaymath}
and
\begin{displaymath}
    [(e,\beta),(e,\beta + \frac{3e}{\delta})]
\end{displaymath}
For the degenerate case $x^2+\delta [x,y]=0$ there are no non-trivial 
extensions outside the diaogonal. 
\par
In the remaining case, where $\delta=0$, the equation is just 
\begin{displaymath}
    x^2+ex=0
\end{displaymath}
and we have a completely different situation. In this case there are 
simple modules of all dimensions. Examples of such are 
\begin{displaymath}
    x\mapsto -ee_{11} \qquad y\mapsto \sum_{i=1}^{n-1}e_{i,i+1}+e_{n,1}
\end{displaymath}
when $e\neq 0$ and
\begin{displaymath}
    x\mapsto e_{n1} \qquad y\mapsto \sum_{i=1}^{n-1}e_{i,i+1}
\end{displaymath}
when $e=0$. 

\subsection{The cusp}\label{cusp}

Consider the ordinary cusp, i.e. the $k$-algebra 
\begin{displaymath}
    A=k\langle x,y\rangle /(y^2-x^3)
\end{displaymath}
A 2-dimensional representation $\phi:A\rightarrow M_{2}(\Gamma)$ maps 
$y^2-x^3$ to 
\begin{displaymath}
    t_{Y}Y-(t_{X}^2-d_{X})X-(d_{Y}-t_{X}d_{X})I
\end{displaymath}
where $X=\phi(x)$ and $Y=\phi(y)$. If the representation is 
simple then $t_{Y}=t_{X}^2-d_{X}=d_{Y}-t_{X}d_{X}=0$, i.e.
$t_{Y}=0$, $d_{X}=t_{X}^2$ and $d_{Y}=t_{X}^3$. The 2-dimensional 
simple modules forms a subvariety of the affine 2-space $Spec 
(k[t_{X},t_{XY}])$, given as the complement of the hypercusp 
$t_{XY}^2-3t_{X}^5=0$, corresponding to the Formanek element. 
\par
Let $V$ be a 1-dimensional simple module, represented by the closed point 
$(a,b)$ on the curve, and let $W$ be another 1-dimensional simple $A$-module. 
Then \newline $Ext_{A}^1(V,W)\neq 0$ if and only if $W$ is 
represented by a closed point $(\omega a,-b)$ on the curve, with 
$\omega$ 
any 
primitive qube root of 1.
\par
Notice that the center $Z(A)\subset A$ is the subalgebra of 
$A$ generated by $t:=x^3=y^2$. It is clear that any surjectiv 
homomorphism of $k$-algebras,
\begin{displaymath}
    \rho_v: A\longrightarrow End_k(V)
\end{displaymath}
will map $Z(A)=k[t]$ into $Z(End_k(V))\simeq k$, inducing a point $v\in 
Simp(k[t])=\mathbf A^1$. Thus $Simp_n(A)$
is fibred over the affine line $Spec(k[t])=\mathbf A^1$. Let 
$\rho_v(x)^3=\rho_v(y)^2=\kappa(v) I_{n}$, where
$\kappa(v)$ is a function on the curve, such that $\kappa(v)=0$ if and 
only if $v=origin=:\underline o$.
Consider the diagram:
\begin{displaymath}
    \diagram
    k[t=x^3=y^2]\dto\drto\\
    A\dto\rto^{\rho_v} &End_k(V) \\
    k[x]/(x^3-\kappa(v))\ast k[y]/(y^2-\kappa(v))\urto
    \enddiagram
\end{displaymath}
Clearly, if $\kappa(v)\neq 0$ the simple representations of $A$ are 
fibered on the cusp with fibres being the simple
representations of $U:=k[x]/(x^3-\kappa(v))\ast 
k[y]/(y^2-\kappa(v))$, isomorphic to the group algebra of
the projective modular group $PSL(2,\mathbf Z)$.
\par
For the 2-dimensional representations as described above, this 
correspond to simple representations
\begin{displaymath}
    x\mapsto\2matr{x_{11}}{x_{12}}0{\omega x_{11}}
    \qquad
    y\mapsto\2matr{y_{11}}0{y_{21}}{-y_{11}}
\end{displaymath}
where $x_{11}^3=y_{11}^2=\kappa(v)$ and $\omega$ is a primitive third 
root of unity. The fibre of $Simp_{2}(A)$ above $v\in{\mathbf A}^1$ 
is an open set of (the three lines) $t_{X}^3-\kappa(v)=0$.
\par
For $\kappa(v)=0$ the matrices,
$\rho_{\underline o}(x)$ and $\rho_{\underline o}(y)$ must both be 
nilpotent $2\times 2$-matrices, and generating $M_{2}(k)$. It turns out 
that we may assume,
$$
\rho_{\underline o}(x)=\left(\begin{array}{cc}
0&\lambda\\
0&0
\end{array}
\right)
,\
\rho_{\underline o}(y)=\left(\begin{array}{cc}
0&0\\
1&0
\end{array}
\right).
$$
where $\lambda\in k$.

A formal versal family is given by,
$$
\tilde\rho(x)=\left(\begin{array}{cc}
t&1+\mu\\
0&t
\end{array}
\right)
,\
\tilde\rho(y)=\left(\begin{array}{cc}
0&0\\
1+\mu&0
\end{array}
\right),
$$
defined over the complete $k$-algebra,
$$
H^A(V)=k[[t,\mu]]/(t^3, (\mu+1)t^2)
$$
For $\mu\neq -1$ this determines a simple representation.

\subsection{Elliptic curve}\label{elliptic curve}

In this section we shall study the geometry of a noncommutative version of an 
elliptic curve given by the $k$-algebra 
\begin{displaymath}
    A=k\langle x,y\rangle/(y^2-x^3-ax-b+q[x,y])
\end{displaymath}

The 1-dimensional simple modules are represented by the 
elliptic curve 
\begin{displaymath}
    C=Simp_{1}(A)=Spec(k[x,y]/(y^2-x^3-ax-b))
\end{displaymath}
A 2-dimensional representation 
\begin{displaymath}
    \phi:A\longrightarrow M_{2}(\Gamma)
\end{displaymath}
is determined by two $2\times 2$-matrices $X=\phi(x)$ and $Y=\phi(y)$, 
satisfying the equation $Y^2-X^3-aX-b+q[X,Y]=0$. Using elementary 
arithmetic for $2\times 2$-matrices, stated in Lemma 
\ref{traceformulas}, it is easily seen that this is 
equivalent to
\begin{align*}
    2qXY+(t_{Y}-qt_{X})Y&+(-t_{X}^2+d_{X}-a-qt_{Y})X\\
    &+(-d_{Y}+t_{X}d_{X}-b-qt_{XY}+qt_{X}t_{Y})=0
\end{align*}
If $\phi$ corresponds to a simple module the set $\{1,X,Y,XY\}$ is 
linearily independent over $k$ and the relations 
\begin{displaymath}
    2q=t_{Y}-qt_{X}=-t_{X}^2+d_{X}-a-qt_{Y}
    =-d_{Y}+t_{X}d_{X}-b-qt_{XY}+qt_{X}t_{Y}=0
\end{displaymath}
holds. Substituting $q=0$ into the rest of the equations gives
\begin{displaymath}
    t_{Y}=-t_{X}^2+d_{X}-a=-d_{Y}+t_{X}d_{X}-b=0
\end{displaymath}
and the trace ring is clearly generated by $t_{X}$ and $t_{XY}$.
The simple locus is given by the non-vanishing of the Formanek element 
\begin{align*}
    d_{[X,Y]}
    &=-t_{XY}^2+(t_{X}^2-4t_{X}^2-4a)(-t_{X}^3-at_{X}+b)\\
    &=-t_{XY}^2+3t_{X}^5+7at_{X}^3-3bt_{X}^2+4a^2t_{X}-4ab
\end{align*}
corresponding to the complement of a hyperelliptic curve of genus 2.
\par
The main tool for studying the indecomposable modules of $A$ is the 
noncommutative Jacobian of the curve. Let $P=(u,v)$ be a 
generic point on the elliptic curve. Then vanishing of the Jacobian matrix 
\begin{displaymath}
    J((f);p)=\matrH {-u^2-ux-x^2-a+qv-qy}{v+y+qx-qu}=(0)
\end{displaymath}
produces a quadratic relation in $x$. There are two different 
solutions outside of the zero set of the discriminant 
\begin{displaymath}
    {\mathcal D}=(u-q^2)^2-4(u^2+q^2u+a-2qv)=0
\end{displaymath}
The two solutions of the quadratic equation can be written in the form 
\begin{displaymath}
    x_{1}=-\frac 12 u+\frac 12 q^2 +\frac 12 \sqrt{\mathcal D}
    \quad
    x_{2}=-\frac 12 u+\frac 12 q^2 -\frac 12 \sqrt{\mathcal D}
\end{displaymath}
In both cases we get 
\begin{displaymath}
    y=-v-q(x-u)
\end{displaymath}
Thus the two solutions of the equation $J((f);p)=(0)$ are given by
\begin{displaymath}
    Q_{i}=(x_{i},y_{i})=(u,-v)+(x_{i}-u)(1,-q)\quad\quad i=1,2
\end{displaymath}

Notice that, with respect to the addition rule on the elliptic 
curve, $-P=(u,-v)$ and it is easily seen that the three points 
$-P, Q_{1}$ and $Q_{2}$ are collinear.

Thus we have the following theorem.
\begin{Thm}\label{collinear}
    Let $P=(u,v)$ be a point on an affine elliptic curve, as given 
    above. Let $(P,Q_{i})\in \R$ for $i=1,2$, then 
    \begin{displaymath}
	P=Q_{1}+Q_{2}
    \end{displaymath}
    Moreover, the slope of the line through the points $Q_{1},Q_{2}$ 
    equals $-q$, and this number is independent of choice of the 
    point $P$ on the elliptic curve. 
\end{Thm}

\begin{proof}
    Follows immedeately from the above discussion.
\end{proof}

Notice the following observation. Suppose there 
exist points $P$ and $Q$ on the elliptic curve with non-vanishing 
extension groups 
\begin{displaymath}
    Ext^1_{A}(P,Q)\neq 0,\qquad Ext^1_{A}(Q,P)\neq 0
\end{displaymath}
According to Theorem \ref{main theorem} this a necessary 
condition for the existence of simple 2-dimensional modules close to 
any indecomposable modules given by the unordered pair $(P,Q)$.
\par
The non-vanishing criterion says that there exsist points $Q_{2}$ and 
$P_{2}$ such that 
\begin{displaymath}
    P=Q+Q_{2}\qquad\qquad Q=P+P_{2}
\end{displaymath}
and such that $Ext^1_{A}(P,Q_{2}), Ext^1_{A}(Q_{1},P_{2})\neq 0$.
According to theorem \ref{collinear} we have $Q=-P+t_{1}(1,-q)$ and 
$P=-Q+t_{2}(1,-q)$ or $-P=Q+t_{2}(1,q)$. Thus we must have 
$t_{1}(1,-q)+t_{2}(1,q)=0$ which is impossible unless $q=0$.
\par
Notice that, by a similar argument, there are no 3-cycles 
$P_{1},P_{2},P_{3}$ such that $Ext_{A}^1(P_{i},P_{i+1})\neq 0$, $i\in 
\mathbf{Z}/(3)$ for the 
elliptic curve. In fact the existence of such a cycle gives 
\begin{align*}
    P_{2}&=-P_{1}+t_{1}(1,-q)\\
    P_{3}&=-P_{2}+t_{2}(1,-q)=P_{1}+t_{1}(i,q)+t_{2}(1,-q)\\
    P_{1}&=-P_{3}+t_{3}(1,-q)\\
\end{align*}
which is easily seen to be impossible.


\end{document}